\documentclass{amsart} % Nicer than default article style: less
\usepackage{amsmath,amsthm}     % Handy math stuff, theorem environments.
\usepackage{amssymb}            % Fancy math symbols.
\usepackage{euscript}           % Nice script font.
\usepackage{enumerate,calc}
\usepackage[matrix,arrow,curve,frame]{xy}    % XY-pic diagram pac

\xymatrixcolsep{1.9pc}                          % Adjust size of diagrams.
\xymatrixrowsep{1.9pc}
\newdir{ >}{{}*!/-5pt/\dir{>}}                  % Make better tailed arrows

\def\longfib{\DOTSB\relbar\joinrel\twoheadrightarrow}

\newtheorem{thm}[subsection]{Theorem}
\newtheorem{defn}[subsection]{Definition}
\newtheorem{prop}[subsection]{Proposition}

\newtheorem{cor}[subsection]{Corollary}
\newtheorem{lemma}[subsection]{Lemma}

\theoremstyle{definition}  % Bold headings and Roman body text.
\newtheorem{example}[subsection]{Example}

\newtheorem{remark}[subsection]{Remark}

  {\end{list}}
  
\newcommand{\dfn}{\textbf} % Make defined words bold.

\newcommand{\mdfn}[1]{\dfn{\mathversion{bold}#1}} % Even make math bold

\newcommand{\tens}              {\otimes}               %tensor
\newcommand{\iso}               {\cong}  
\newcommand{\cat}{\EuScript}    % Use \EuScript to name a category.
\newcommand{\cA}{{\cat A}}      % Only seems to work for caps, and only gets
\newcommand{\cB}{{\cat B}}      % first letter.
\newcommand{\cC}{{\cat C}}
\newcommand{\cD}{{\cat D}}

\newcommand{\cI}{{\cat I}}

\newcommand{\cL}{{\cat L}}
\newcommand{\cM}{{\cat M}}

\newcommand{\Top}{{\cat Top}}

\newcommand{\Spectra}{{\cat Spectra}}
\newcommand{\Set}{{\cat Set}}

\newcommand{\sSet}{s{\cat Set}}
\newcommand{\Sch}{{\cat Sch}}

\newcommand{\sPre}{sPre}
\newcommand{\Pre}{Pre}

\newcommand{\Ho}{\text{Ho}\,}

\newcommand{\field}[1]  {\mathbb #1} % Use blackboard bold for these sets
\newcommand{\A}         {\field A}

\newcommand{\Z}         {\field Z}

\DeclareMathOperator*{\colim}{colim}
\DeclareMathOperator*{\hocolim}{hocolim}

\DeclareMathOperator*{\holim}{holim}

\DeclareMathOperator{\spec}{Spec}
\DeclareMathOperator{\Hom}{Hom}
\DeclareMathOperator{\dgn}{dgn}
\DeclareMathOperator{\Tot}{Tot}
\DeclareMathOperator{\homp}{hom_+}

\DeclareMathOperator{\Map}{Map}
\DeclareMathOperator{\hMap}{hMap}

\DeclareMathOperator{\sk}{sk}
\DeclareMathOperator{\cosk}{cosk}
\DeclareMathOperator{\coskp}{cosk}

\DeclareMathOperator{\eq}{eq}

\DeclareMathOperator{\Ex}{Ex}
\DeclareMathOperator{\hgt}{ht\,}

\newcommand{\ra}{\rightarrow}                   % right arrow
\newcommand{\lra}{\longrightarrow}              % long right arrow
\newcommand{\la}{\leftarrow}                    % left arrow
\newcommand{\lla}{\longleftarrow}               % long left arrow
\newcommand{\llra}[1]{\stackrel{#1}{\lra}}      % labeled long right
\newcommand{\we}{\llra{\sim}}                   % weak equivalence

\newcommand{\cof}{\rightarrowtail}              % cofibration
\newcommand{\fib}{\twoheadrightarrow}           % fibration
\newcommand{\trfib}{\stackrel{\sim}{\longfib}}
\newcommand{\trcof}{\stackrel{\sim}{\cof}}

\newcommand{\inc}{\hookrightarrow}              % inclusion
\newcommand{\dbra}{\rightrightarrows}           % double arrow for eqlizer

\newcommand{\blank}{-}                          % A hyphen, as in
\newcommand{\sing}{Sing}
\newcommand{\Sing}{\sing\,}

\newcommand{\ovcat}{\downarrow}

\newcommand{\restr}[1]{\!\mid_{#1}}

\newcommand{\bd}[1]{\partial\Delta^{#1}}

\newcommand{\del}[1]{\Delta^{#1}}

\newcommand{\he}{\simeq}

\newcommand{\Sm}{Sm}

\newcommand{\rea}[1]{|{#1}|}             %geometric realization of #1
\newcommand{\map}{\rightarrow}

\newcommand{\ceck}[1]{\Cech(#1)}         %Cech complex for #1
\newcommand{\oceck}[1]{\Cech^{o}(#1)}    %Ordered Cech complex for #1
\newcommand{\oreal}[1]{\rea{\oceck{U}}}  %Realization of ordered Cech cplex
\newcommand{\creal}[1]{\rea{\ceck{U}}}   %Realization of the Cech complex
\newcommand{\Cech}{\check{C}}
\newcommand{\CCech}{\v{C}ech\ }

\newcommand{\fibr}[1]{{\field H}#1}
\newcommand{\piH}{\pi HC}

\numberwithin{equation}{subsection}

\newenvironment{myequation}
  {\addtocounter{subsection}{1}\begin{eqnarray}}
  {\end{eqnarray}$\!\!$}

\newcommand{\ch}{\check{\cC}}

\begin{document}

\title{Hypercovers and simplicial presheaves}

\author{Daniel Dugger} 
\author{Sharon Hollander}
\author{Daniel C. Isaksen}

\address{Department of Mathematics\\ University of Oregon\\ Eugene, OR
97403}

\address{Department of Mathematics\\ University of Chicago\\ Chicago,
IL 60637}

\address{Department of Mathematics\\ University of Notre Dame\\
Notre Dame, IN 46556}

\email{ddugger@math.uoregon.edu}

\email{sjh@math.uchicago.edu}

\email{isaksen.1@nd.edu}

\begin{abstract}
We use hypercovers to study the homotopy theory of simplicial
presheaves.  The main result says that model structures for simplicial
presheaves involving local weak equivalences can be constructed by
localizing at the hypercovers.  One consequence is that the fibrant
objects can be explicitly described in terms of a hypercover descent
condition, and the fibrations can be described by a relative descent
condition.  
%These ideas are central to constructing realization
%functors on the homotopy theory of schemes \cite{DI1,Is}.  
%We give a few other applications for this new description of the 
%homotopy theory of simplicial presheaves.
We give a few applications for this new description of the 
homotopy theory of simplicial presheaves.
\end{abstract}

\maketitle

\tableofcontents

\section{Introduction}

This paper is concerned with the subject of homotopical sheaf theory,
as it has developed over time in the articles
\cite{I,B,BG,Th,Jo,J1,J2,J3,J4}.  Given a fixed Grothendieck site
$\cC$, one wants to consider contravariant functors $F$ defined on
$\cC$ whose values have a homotopy type associated to them.  The most
basic question is: what should it mean for $F$ to be a sheaf?  The
desire is for some kind of local-to-global property---also called a
{\em descent} property---where the value of $F$ on an object $X$ can
be recovered by homotopical methods from the values on a cover.
Perhaps the earliest instance where such a concept had to be tackled
was in algebraic geometry, where people had to deal with presheaves of
chain complexes defined on a space $X$.  Because of its abelian nature
this could be handled by classical homological algebra, and led to the
Grothendieck definition of hypercohomology.  Much later, people
encountered the non-abelian example of algebraic $K$-theory.  Here the
site $\cC$ is a category of schemes, and the functor $F$ assigns to
each scheme $X$ its algebraic $K$-theory spectrum $K(X)$.  Thomason's
paper \cite{Th} (building on earlier work from \cite{B,BG}) combined
homotopy theory and sheaf theory to study the descent properties of
this functor.

The work of \cite{BG,Jo,J2} brought the use of model categories into
this picture.  In the most modern of these \cite{J2}, Jardine defined
a model category structure on presheaves of simplicial sets with the
property that the weak equivalences are {\em local} in nature.
Classical invariants such as sheaf cohomology arise in this setting as
homotopy classes of maps into certain Eilenberg-Mac Lane objects, and
the whole theory can in some sense be regarded as the study of
non-additive sheaf cohomology.  Jardine's model structure has recently
served as the foundation from which Morel and Voevodsky built their
$\A^1$-homotopy theory for schemes \cite{MV}.

One important ingredient missing from Jardine's work is a description
of the fibrant objects and the fibrations.  They can be characterized
in terms of a certain lifting property, but this is not so
enlightening and not very useful in practice.  Our main goal in this
paper is to give explicit, simple characterizations of the fibrations
and fibrant objects in terms of descent conditions.  This is
equivalent to describing Jardine's model category as a very explicit
Bousfield localization.

\medskip

To explain the basic ideas, let's assume our site is the category of
topological spaces equipped with the usual open covers.  A presheaf of
sets $F$ is a sheaf if $F(X)$ is the equalizer of $\prod_a F(U_a)
\dbra \prod_{a,b} F(U_a \cap U_b)$ whenever $\{U_a\}$ is an open cover
of $X$.  This equalizer is in fact the same as the inverse limit of
the entire cosimplicial diagram
\[ \xymatrix{
   \prod_a F(U_a) \ar@<0.5ex>[r]\ar@<-0.5ex>[r] &
    \prod_{a,b} F(U_{ab}) \ar@<0.6ex>[r]\ar[r]\ar@<-0.6ex>[r] &
    \cdots
}
\]
where we have abbreviated $U_{a_0\ldots a_n}$ for $U_{a_0}\cap \cdots
\cap U_{a_n}$ and have refrained from drawing the codegeneracies for
typographical reasons.  For a presheaf of {\it simplicial\/} sets (or
taking values in some other homotopical objects like {\it spectra\/}),
it is natural to replace the limit by a homotopy limit.  So one
requires that $F(X)$ be weakly equivalent to the homotopy limit of the
above cosimplicial diagram.  This property, when it holds for all open
covers, is called \dfn{\CCech descent}.  It can also be expressed in a
slightly more compact way, if one recalls that the \CCech complex
$\Cech{U}$ associated to a cover $\{U_a\}$ of $X$ is the simplicial
object $[n] \mapsto \coprod_{a_0 \cdots a_n} U_{a_0\cdots a_n}$.  Then
$F$ satisfies \CCech descent if the natural map
\[ F(X) \ra \holim_n F(\Cech{U}_n) 
\]
is a weak equivalence.

A motivating example is given by the functor $\Top^{op} \ra \Spectra$
taking $X$ to $E^X$, where $E$ is a fixed spectrum and $E^X$ denotes
the function spectrum.  This functor has \CCech descent, because $X$
is weakly equivalent to the homotopy colimit of the \CCech complex 
associated to 
any open cover of $X$---this was shown in \cite[Thm. 1.1]{DI1}.

Now, it is not true that the fibrant objects in Jardine's model
category are just the simplicial presheaves which satisfy \CCech
descent (although this erroneous claim has appeared in a couple of
preprints, for instance \cite{HS}).  See the appendix, Example
\ref{ex:counter}, for an example.  What we show in this paper is that
one has to instead consider descent for all {\it hypercovers\/}.  A
hypercover is a simplicial object $U$, augmented by $X$, which is
similar to a \CCech complex except in level $n$ we only need to have a
{\it cover} of the $n$-fold intersections $U_{a_0\cdots a_n}$.  A
precise definition requires a morass of machinery (see Section
\ref{se:hc}).  A simplicial presheaf $F$ {\it satisfies descent\/} for
the hypercover $U \ra X$ if the natural map
\[ F(X) \ra \holim_n F(U_n) 
\]
is a weak equivalence (see Definition~\ref{de:hc-desc}).

What we will show is that the fibrant objects in Jardine's model
category are essentially the simplicial presheaves which satisfy
descent for all hypercovers (there is an analagous criterion for
fibrations in terms of a relative descent property, given in
Section~\ref{se:fibs}):

\begin{thm}
\label{th:J-fib}
The fibrant objects in Jardine's model category $\sPre(\cC)_\cL$ are
those simplicial presheaves that:
\begin{enumerate}[(1)]
\item are fibrant in the injective model structure $\sPre(\cC)$, and
\item satisfy descent for all hypercovers $U \ra X$.
\end{enumerate}
\end{thm}

The {\em injective} model structure on $\sPre(\cC)$ just refers to
Jardine's model structure for the discrete topology on $\cC$
(see Section~\ref{se:mc} for more about this).  The fibrancy
conditions for this model structure are awkward to describe, but they
also aren't very interesting---they have no dependence on the
Grothendieck topology, only on the shape of the underlying category
$\cC$.  The conditions require that each $F(X)$ be a fibrant
simplicial set, certain maps $F(X)\ra F(Y)$ be fibrations, and more
complicated conditions of a similar `diagrammatic nature'.  In
practice such conditions are not very important, and in fact there's a
way to get around them completely by using the {\em projective}
version of Jardine's model structure; see Theorem \ref{th:U-fib} below
and the discussion in Section~\ref{se:mc}.

We'd like to point out that the above theorem can be re-interpreted in
terms of giving `generators' and `relations' for the homotopy theory
of simplicial presheaves, in the manner introduced by \cite{D}.  Using
the language of that paper, we prove

\begin{thm} 
Jardine's model category $\sPre(\cC)_{\cL}$ is Quillen
equivalent to the {\em universal homotopy theory} $U\cC/S$ constructed by
\begin{enumerate}[(1)]
\item formally adding homotopy colimits to the category $\cC$, to
create $U\cC$;
\item imposing relations requiring that for every hypercover
$U \ra X$, the map $\hocolim_n U_n \map X$ is a weak equivalence.
\end{enumerate}
\end{thm}

In other words, the result says that everything special about the
homotopy theory of simplicial presheaves can be derived from the basic
fact that one can reconstruct $X$ as the homotopy colimit of any of
its hypercovers.  The above theorem is crucial to the construction of
\'etale realization functors for $\A^1$-homotopy theory \cite{Is}, as
well as the analogous question about topological realization functors
\cite{DI1}.

One advantage of the model structure $U\cC/S$ over the model structure
$\sPre(\cC)_\cL$ is that the fibrant objects are much easier to
describe.  The inexplicit fibrancy conditions of the injective model
structure are replaced by a much simpler condition.  Compare the
following result to Theorem~\ref{th:J-fib}.

\begin{thm}
\label{th:U-fib}
The fibrant objects in the model category $U\cC/S$ are
those simplicial presheaves that:
\begin{enumerate}[(1)]
\item are objectwise-fibrant ({\em i.e.}, each $F(X)$ is a fibrant
simplicial set),
\item satisfy descent for all hypercovers $U \ra X$.
\end{enumerate}
\end{thm}

Again, an analagous criterion for all fibrations is given in
Section~\ref{se:fibs}.

The main ideas we use to prove these results are very simple, and
worth summarizing. They exactly parallel classical facts about
CW-complexes.  The two key ingredients are:
\begin{enumerate}[(i)]
\item In the category of simplicial presheaves one can construct
objects analogous to CW-complexes, the only difference being that one
has different kinds of $0$-simplices corresponding to the different
representable presheaves $rX$.  (And as a consequence, there are
different kinds of $n$-simplices corresponding to the objects
$\del{n}\tens rX$.)  Every simplicial presheaf has a cellular
approximation built up out of representables in this way (see
\cite[Section 2.6]{D}).
\item Weak equivalences for simplicial presheaves are characterized by
a certain `local lifting criterion', where lifting problems can
be solved by passing from a representable object to the pieces of a
cover.  See Proposition~\ref{pr:loclift} and the paper \cite{DI2}.
\end{enumerate}
From these two basic principles, it's inevitable that hypercovers will
arise in the solution of lifting problems.  One starts
building a lift inductively on a CW-approximation, and the
obstructions to extending the lift are made to vanish by passing to a
finer cover at each stage.  Thus, one finds oneself inductively
constructing a hypercover.  These ideas are explored in detail in
Section~\ref{se:hclift}.

\medskip

This paper came into existence because we needed to use
Theorems~\ref{th:hcloc}, \ref{th:Verdier}, and \ref{th:allagree}(c,d)
in other work.  We originally hoped for a very short paper, but to
actually write down complete proofs one has to be able to manipulate
hypercovers with a certain amount of ease---and the literature on this
subject is not the most helpful.  So in the end a large portion of the
paper has been devoted to carefully setting down the machinery of
hypercovers, hopefully in a way that will be usable by other people.
For this reason the paper sometimes takes on an expository tone.  We
have tried to be clear and thorough, and for good or bad this has come
at the expense of brevity.  Also, one of our goals has been to adopt
definitions which can be applied to any Grothendieck site, not just
the classical ones which get used most often.  The result is theorems
which are simple enough to state and prove, but sometimes hard to
apply in practice.  To complement this, we have included the
reductions to Verdier sites (Section \ref{se:Verdier}) and internal
hypercovers (Section \ref{se:internal}) that one can implement for sites
like those encountered in algebraic geometry.  This subject of
homotopical sheaf theory is rapidly finding applications in many
contexts, so we have tried to give a presentation that is clear
enough, and general enough, to be useful to a variety of
practitioners.

\subsection{Organization of the paper}\mbox{}\par
In Section \ref{se:mc} we review the basic model categories that will
be used throughout the paper.  One of these is Jardine's model
structure, and the other is a Quillen equivalent version which has
fewer cofibrations and more fibrations.  We assume throughout that the
reader is familiar with the theory of model categories---the original
reference for this subject is \cite{Q}, but we generally follow
\cite{H} in notation and terminology. \cite{Ho} is also a good
reference.  

Section \ref{se:we} reviews material from \cite{DI2} on lifting
properties for simplicial presheaves, and how these can be used to
characterize local weak equivalences.  Section \ref{se:hc} introduces
the machinery needed for defining and working with hypercovers.  The
section is a bit long, and serves mostly as a reference section for
the rest of the paper---it can comfortably be skimmed the first time
through.

In Section \ref{se:hclift} we show how hypercovers enter into the
solution of lifting problems in the homotopy theory of simplicial
presheaves.  These are the key observations which are needed for the
main results.  The proofs of these main results are then given in
Sections~\ref{se:hcloc} and \ref{se:fibs}, where they appear as 
Theorem~\ref{th:hcloc}, Corollary~\ref{co:hcloc}, and
Theorem~\ref{th:fibrations}.

One application of the results on hypercovers is to realization
functors from the homotopy theory of schemes---this is treated in the
papers \cite{DI1,Is}.  In Section~\ref{se:app} we give a few more
applications.  One of the most interesting, given in
Section~\ref{subse:changesite}, is a much simpler approach to the
change-of-site functors of \cite{MV}.  We also discuss a
generalization of the Verdier hypercovering theorem in
Theorem~\ref{th:Verdier}.

In applications one rarely wants to work with {\it all\/} hypercovers,
because this is just too broad a class of objects.  In the case of the
`geometric' sites which are most commonly used, one can adopt more
restrictive definitions and have all the above results still go
through.  These reductions are explored in Sections~\ref{se:Verdier}
and~\ref{se:internal}.  We axiomatize what is necessary into the
notion of a {\it Verdier site\/}, which comes equipped with a special
class of `basal hypercovers'.  These ideas appear sporadically in
Section~\ref{se:app}, but hopefully the reader can just refer back to
the later sections as necessary.

Finally, the paper contains an appendix which explores the difference
between \CCech descent and hypercover descent.  Again, the principal
motivation comes from the fact that \CCech descent is more easily
dealt with in practice.  We show, among other things, that having
descent for \CCech complexes is equivalent to having descent for all
{\it bounded} hypercovers (the ones where the refinement process stops
at some finite level).  This is an important ingredient in \cite{DI1}.

\subsection{Notation and Terminology}

If $X$ is an object of a site $\cC$, then the \dfn{representable}
simplicial presheaf $rX$ on $\cC$ is given by the formula $rX(Y) =
\Hom_{\cC}(Y, X)$.  Note that each simplicial set $rX(Y)$ is discrete.
If $U$ is a simplicial object of $\cC$, then $rU$ is the simplicial
presheaf given by the formula $rU(Y)_n = \Hom_{\cC} (Y, U_n)$---these,
of course, are usually not discrete.  We frequently abuse notation and
write simply $X$ (or $U$) for the presheaf $rX$ (or $rU$).  

If $S$ is a scheme, then $\Sch/S$ denotes the category of schemes of
finite-type over $S$.  The full subcategory of schemes which are
smooth over $S$ is denoted $\Sm/S$.  Finally, in a simplicial model
category we write $\Map(A,B)$ for the simplicial mapping space.

%%%%%%%%%%%%%%%%%%%%%%%%%%%%%%%%%%%%%%%%%%%%%%%%%%%%%%%%%%%%%%%%%%%%%%

\section{Model structures on simplicial presheaves}
\label{se:mc}

We start by recalling that for any small category $\cC$ there are two
Quillen equivalent model structures on the category of diagrams
$\sSet^\cC$.  In each case a map $D \ra E$ is a weak equivalence if
$D(c) \ra E(c)$ is a weak equivalence of simplicial sets for each $c$ in
$\cC$.  Such a map is usually called an \dfn{objectwise weak
equivalence}.  In the \dfn{projective model structure} on $\sSet^\cC$
one defines a map $D\ra E$ to be
\begin{enumerate}[(1)]
\item A {\em fibration} if every $D(c) \ra E(c)$ is a fibration of simplicial
sets ({\em i.e.}, $D\ra E$ is an objectwise fibration).
\item A {\em cofibration} if it has the left-lifting-property with respect
to the acyclic fibrations.
\end{enumerate}
Dually, in the \dfn{injective model structure} the cofibrations are
objectwise and the fibrations have the right-lifting-property with
respect to acyclic cofibrations.  The names `projective' and
`injective' come from the analogy between the two usual model
structures on chain complexes of $R$-modules.  For notational
convenience, the projective model structure is denoted $U\cC$ (as was
done in \cite{D}, where it was pointed out that $U\cC$ has a certain
universal property) and the injective model structure is denoted
$\sPre(\cC)$.

When $\cC$ comes equipped with a Grothendieck topology, then one can
construct refinements of these model structures which reflect the
topology on $\cC$.  A map of simplicial presheaves $F\ra G$ is a
\dfn{local weak equivalence} if it induces isomorphisms on all sheaves
of homotopy groups~\cite{I,Jo,J2}.  In this paper we will use an
alternative characterization in terms of homotopy liftings, described
below.

Jardine's model structure on simplicial presheaves is the left
Bousfield localization of $\sPre(\cC)$ at the class $\cL$ of local
weak equivalences; we denote this localization as
$\sPre(\cC)_{\cL}$.  Of course since $\cL$ is a {\it class} of maps
there is no {\it a priori\/} guarantee that the Bousfield localization
exists, but Jardine was able to construct it directly---it is only
after the fact that one can identify it as a localization.

Similarly, one can define a model structure $U\cC_{\cL}$ by localizing
$U\cC$ at the same class $\cL$ (cf. \cite[Thm.~1.5]{Bl}).  The identity
maps induce a Quillen equivalence $U\cC_{\cL} \ra \sPre(\cC)_{\cL}$, so
once again these are projective and injective versions of the same
underlying homotopy theory.  The injective version has the advantage
that every object is cofibrant, but in the projective version the
fibrant objects are easier to understand and the representable
presheaves are still cofibrant.  Also, it is usually easier to construct
functors out of the projective version \cite{D}.
We state most of our results only in terms of $\sPre(\cC)_{\cL}$,
but analogous statements for $U\cC_{\cL}$
are also true with only 
minor differences between the proofs.  

Both $U\cC$ and $\sPre(\cC)$ are proper, simplicial model
categories: if $F$ is a simplicial presheaf and $K$ is a simplicial
set then $K\tens F$ and $F^K$ are defined objectwise, by
\[ (K\tens F)(X) = K \times F(X) \qquad\text{and}\qquad
   (F^K)(X)=F(X)^K.
\]
From general considerations \cite[Thm. 4.1.1]{H}, all localizations of
$U\cC$ and $\sPre(\cC)$ that we consider are also left proper,
simplicial model categories.

\begin{remark}
\label{re:hocolim}
If $F$ is a simplicial presheaf, then one obtains a diagram $D_F
\colon\Delta^{op} \ra \sPre(\cC)$ by sending $[n]$ to $F_n$.  Here
$F_n$ is just a presheaf of sets, but we can regard it as a discrete
simplicial presheaf in the obvious way.  The realization of this
simplicial diagram is precisely $F$.  The Bousfield-Kan map $\hocolim
D_F \ra \rea{D_F}$ is a weak equivalence in this case, by some basic
model category theory.  So any simplicial presheaf $F$ is weakly
equivalent to $\hocolim D_F$. This observation will be needed often.
\end{remark}

\section{Local weak equivalences and local lifting properties}
\label{se:we}

Local weak equivalences are usually defined in terms of sheaves of
homotopy groups.  Here we recall a different description which is more
suitable for our purposes.  See \cite{DI2} for the proof that the two
definitions agree and for more details on the results in this section.

First, recall that if $X$ is in $\cC$ and $F$ and $G$ are simplicial
presheaves, then a diagram such as
\[ \xymatrix{
\Lambda^{n,k} \tens X \ar[r] \ar[d] & F \ar[d] \\
\del{n} \tens X \ar[r] & G}
\]
has \dfn{local liftings} if there exists a covering sieve $R$ of $X$
such that for every map $U\ra X$ in the sieve, the diagram one obtains
by restricting from $X$ to $U$ has a lifting $\del{n}\tens U \ra F$.
These liftings are not required to be compatible for the
different $U$'s.  A map $F\ra G$ is called a \dfn{local fibration} if
it has local liftings with respect to the maps $\Lambda^{n,k}\tens X
\ra \del{n}\tens X$, for all $X$ in $\cC$.  A simplicial presheaf is
called \dfn{locally fibrant} if $F\ra *$ is a local fibration.

\begin{prop}[{\cite[Th. 6.15]{DI2}}]
\label{pr:loclift}
A map $F\ra G$ between locally fibrant simplicial presheaves is a
local weak equivalence if and only if every square
\[ \xymatrix{
\bd{n}\tens X \ar[r] \ar[d] & F \ar[d] \\
\del{n} \tens X \ar[r] & G}
\]
has local \dfn{relative homotopy-liftings}, in the following sense: after
restricting to the pieces $U\ra X$ of some covering sieve, one has
maps $\del{n}\tens U \ra F$ making the upper triangle commute on the
nose and the lower triangle commute up to simplicial homotopy relative
to $\bd{n}\tens U$.
\end{prop}

The reader may consult \cite{DI2} for a detailed discussion of this
kind of relative-homotopy-lifting property.

The following two results from \cite{DI2} will be used later.  Recall
that a map is a \dfn{local acyclic fibration} if it is both a local
fibration and a local weak equivalence.

\begin{prop}[{\cite[Prop. 7.2]{DI2}}]
\label{pr:trlfib}
A map $F \map G$ admits local liftings in 
every square
\[ \xymatrix{ \bd{n} \tens X \ar[r]\ar[d] & F \ar[d] \\
              \del{n} \tens X \ar@{.>}[ur]\ar[r] & G
}\]
if and only if it is a local acyclic fibration.
\end{prop}

One consequence of the above result is that local acyclic fibrations
are closed under pullbacks (in \cite{J2} this was proven only when the
domain and codomain are locally fibrant).

\begin{prop}[{\cite[Lemma 19]{J4}},{\cite[Cor. 7.4]{DI2}}]
\label{pr:fsm7}
Let $F\ra G$ be a local fibration (resp. local acyclic fibration).  
If $K\inc L$ is an
inclusion of finite simplicial sets, then the induced map
\[    F^L \ra F^K \times_{G^K} G^L 
\]
is a local fibration (resp. local acyclic fibration).
\end{prop}

Let $f\colon E \ra B$ be a map between presheaves of sets. One says
that $f$ is a \dfn{generalized cover} (or {\it local epimorphism\/})
if it has the following property: given any map $rX \ra B$, there is a
covering sieve $R\inc X$ such that for every element $U\ra X$ in $R$,
the composite $rU \ra rX \ra B$ lifts through $f$.  The `generalized'
adjective is there to remind us that we are looking at a map between
presheaves, not actual objects of the site.  In the case where $B$ is
representable and $E$ is a coproduct $\coprod E_a$ of representables,
$f$ is a generalized cover precisely when the sieve generated by
the maps $\{E_a \ra B\}$ is a covering sieve of $B$.

For a simplicial presheaf $F$, let $\tilde{M}_nF$ denote the $0$th
object of $F^{\bd{n}}$ (the `tilde' is to distinguish this from a
slightly different construction used later in the paper).  This is the
presheaf of sets whose value $\tilde{M}_nF(X)$ is the set of all maps
$\bd{n}\ra F(X)$.  There is a natural map $F_n\ra \tilde{M}_n F$
induced by $F^{\del{n}} \ra F^{\bd{n}}$.  Proposition~\ref{pr:trlfib}
can be rephrased as saying that $F\ra G$ is a local acyclic fibration
if and only if the maps
\begin{myequation}
\label{eq:matchingmap}
 F_n \ra \tilde{M}_n F \times_{\tilde{M}_n G} G_n 
\end{myequation}
are generalized covers, for all $n\geq 0$.  Using this observation,
most properties of generalized covers can automatically be seen to
hold for local acyclic fibrations.

%%%%%%%%%%%%%%%%%%%%%%%%%%%%%%%%%%%%%%%%%%%%%%%%%%%%%%%%%%%%%%%%%%%%

\section{Background on hypercovers}
\label{se:hc}

This section contains the necessary machinery for defining and working
with hypercovers.  Unfortunately there is quite a bit of annoying
category theory, and some readers may wish to only skim this section
their first time through.  This should be enough to understand
the basic definitions that are used throughout the paper.  In Section
\ref{se:cosk} we recall the coskeleton and degeneration functors,
which appear when passing between simplicial objects and truncated
simplicial objects.  These notions are used later in the paper, but
only in fairly technical contexts.

\subsection{The definition}

\begin{defn}
\label{de:hypercover}
Let $X$ belong to $\cC$ and suppose that $U$ is a simplicial presheaf
with an augmentation $U \ra X$.  This map is called a \dfn{hypercover}
of $X$ if each $U_n$ is a coproduct of representables, and $U\ra X$ is
a local acyclic fibration.
\end{defn}

Using (\ref{eq:matchingmap}) one can rewrite the second condition in a
more explicit way: it says that the maps $U_0\ra X$, $U_1\ra
U_0\times_X U_0$, and $U_n \ra \tilde{M}_n U$ (for $n\geq 1$) are all
generalized covers.  This is not particularly enlightening, but it's
easy to provide some intuition behind it. For convenience we assume our
Grothendieck topology is given by a basis of covering families. Then
the easiest examples of hypercovers are the \CCech complexes, which
have the form
\[ \xymatrix{  
\cdots \coprod U_{a_0 a_1 a_2} \ar[r]<0.6ex>\ar[r]\ar[r]<-0.6ex> 
    &\coprod U_{a_0 a_1} \ar[r]<0.5ex>\ar[r]<-0.5ex>
    &\coprod U_{a_0} \ar[r]
    & X
}
\]
for some chosen covering family $\{U_a \ra X\}$.  Here $U_{a_0 \cdots
a_n}$ is the fibre-product $U_{a_0} \times_X \cdots \times_X U_{a_n}$.
The \CCech complexes are the hypercovers for which the maps $U_1\ra
U_0\times_X U_0$ and $U_n \ra \tilde{M}_nU$ are all {\it
isomorphisms}.  In an arbitrary hypercover one takes the iterated
fibre-products at each level but then is allowed to refine that object
further, by taking a generalized cover of it.  We refer the reader to
\cite[Section 8]{AM} for further discussion of hypercovers.

Next is the formal definition of hypercover descent:

\begin{defn}
\label{de:hc-desc}
An objectwise-fibrant simplicial presheaf $F$ \dfn{satisfies descent}
for a hypercover $U \map X$ if the natural map from $F(X)$ 
to the
homotopy limit of the diagram
\[ \xymatrix{
   \prod_a F(U^a_0) \ar@<0.5ex>[r]\ar@<-0.5ex>[r] &
    \prod_{a} F(U^a_{1}) \ar@<0.6ex>[r]\ar[r]\ar@<-0.6ex>[r] &
    \cdots
}
\]
is a weak equivalence.  Here the products range over the representable
summands of each $U_n$.  If $F$ is not objectwise-fibrant, we say it
satisfies descent if some objectwise-fibrant replacement for $F$ does.
\end{defn}

The definition has been arranged so that if $F\ra G$ is an objectwise
weak equivalence, then $F$ satisfies descent for $U\ra X$ if and only
if $G$ does.  While the definition reflects our intuitive notion of
descent, the next lemma gives a more concise reformulation in terms of
simplicial mapping spaces.

\begin{lemma}\mbox{}\par
\label{le:descent}
\begin{enumerate}[(i)]
\item
A simplicial presheaf $F$ satisfies descent for a hypercover $U\ra X$
if and only if $\Map(X, \hat{F}) \ra \Map(U, \hat{F})$ is a weak
equivalence of simplicial sets, where $\hat{F}$ is an injective-fibrant
replacement for $F$.
\item Let
$U'$ be a cofibrant replacement for $U$ in $U\cC$.  Then $F$
satisfies descent for $U\ra X$ if and only if $\Map(X, \hat{F}) \ra
\Map(U', \hat{F})$ is a weak equivalence of simplicial sets, where
$\hat{F}$ is an objectwise-fibrant replacement for $F$.
\end{enumerate}
\end{lemma}

Note that any split hypercover (see Definition~\ref{de:split}) is
cofibrant in $U\cC$, in which case one can apply (ii) with
$U'=U$.

\begin{proof}
This is by general nonsense.  Consider the diagram $\Delta^{op} \ra
\sPre(\cC)$ given by $[n] \ra U_n$, and let $\tilde{U}$ be its
homotopy colimit.  This is not the same as $U$, but there is a map
$\tilde{U} \ra U$ which is an objectwise weak equivalence (see
Remark~\ref{re:hocolim}).  Let $\hat{F}$ be an injective-fibrant
replacement for $F$, which {\it a fortiori\/} is an objectwise-fibrant
replacement as well.  Then $\Map(\tilde{U}, \hat{F})$ is weakly
equivalent to $\Map(U, \hat{F})$ since $\tilde{U} \ra U$ is a weak
equivalence between injective-cofibrant objects.  But $\Map(\tilde{U},
\hat{F})$ is
\[ \Map(\hocolim_n U_n,\hat{F}) \he \holim_n \Map(U_n,\hat{F}) \he
\holim_n {\textstyle{\prod_a}} \hat{F}(U_n^a).
\]
Since $\Map(X, \hat{F})$ is equal to $\hat{F}(X)$, 
the condition that
$\Map(X,\hat{F})\ra \Map(U,\hat{F})$ be a weak equivalence is a direct
translation of the homotopy limit formulation in
Definition~\ref{de:hc-desc}.  This proves (i).

For (ii), note that each $U_n$ is cofibrant in $U\cC$ and so
$\tilde{U}=\hocolim_n U_n$ is also cofibrant.  In other words
$\tilde{U}$ is a cofibrant replacement for $U$, and so $\tilde{U}\he
U'$.  If $\hat{F}$ is an objectwise replacement for $F$ then
$\Map(U',\hat{F})\he \Map(\tilde{U},\hat{F})$, and as in (i) the
latter is equivalent to $\holim_n \prod_a \hat{F}(U_n^a)$.  The rest
of the proof is the same.
\end{proof}

A more elegant way to phrase the above result is to say that $F$
satisfies descent for $U\ra X$ if and only if $\hMap(X,F) \ra
\hMap(U,F)$ is a weak equivalence of simplicial sets, where
$\hMap(\blank,\blank)$ denotes a homotopy function complex
\cite[Ch.~17]{H} in either $U\cC$ or $\sPre(\cC)$.

\subsection{Machinery}\mbox{}\par

Definition~\ref{de:hypercover} is very compact, but it's not always
such an easy thing to work with.  For the rest of this section we will
set down more convenient techniques for constructing and working with
hypercovers.  This material is used throughout the paper, but many
readers will want to skip ahead and refer back to this section only
when needed.

\medskip

Let $\cM$ be a category which is complete and co-complete---in our
applications $\cM$ is $\Pre(\cC)$, but for the moment let us work
in the more general setting.  Let $\Delta_+$ denote the augmented
cosimplicial indexing category: it is obtained by adjoining an initial
object $[-1]$ to $\Delta$.  Let $s_+\cM$ denote the category of
functors $\Delta^{op}_+ \ra \cM$, i.e. the category of augmented
simplicial objects.  We regard a simplicial set $K$ as belonging
to $s_+\Set$ by letting $K_{-1}$ consist of a single point.

If $S$ is a set and $X$ belongs to $\cM$, let $X^{S}$ denote a
product of copies of $X$ indexed by the elements of $S$.  Given a
simplicial set $K$ and an object $W$ of $s_+\cM$, we regard these as
functors $K\colon \Delta^{op}_+ \ra \Set$ and $W\colon\Delta^{op}_+
\ra \cM$ and then form the resulting end, denoted $\homp(K,W)$:
\[ \homp(K,W):=\eq \Bigl [ \prod_n W_n^{K_n} \dbra 
                           \prod_{[n]\ra [m]} W_m^{K_n} \Bigr ].
\]
The $+$ subscript is to remind us of the augmentations.

\begin{remark}
\label{re:adjoint}
As with any end, this construction exhibits a useful adjointness
property.  If $Z$ is in $\cM$, then the maps $Z \ra \homp(K,W)$ in $\cM$
correspond bijectively with the maps $Z\tens K \ra W$ in $s_+\cM$.
Here $Z\tens K$ is the augmented simplicial object which in dimension
$n$ is a coproduct, indexed by the set $K_n$, of copies of $Z$.
\end{remark}

In the unaugmented simplicial category $s\cM$, we can compute
unaugmented ends $\hom(K, W)$ in an analogous way.  Again, this
construction is right adjoint to tensoring with $K$.

The following lemma can be proved with the above adjointness
property and the Yoneda lemma.

\begin{lemma}
\label{le:hompo}
Let $W \ra X$ be an augmented simplicial object (that is, $X$ is the
augmentation).  
\begin{enumerate}[(i)]
\item $\homp(K,W) \iso \hom^X(K,W)$,
where $\hom^X(K, W)$ is computed in the unaugmented simplicial overcategory
$s(\cM\ovcat X)$.
\item $\homp(K,W)\iso \hom(K,W)$ if $K$ is connected.
\item $\homp(\emptyset,W)\iso X$, and so for any simplicial set $K$
there is a canonical map $\homp(K,W) \ra X$.
\item $\homp(\del{n},W)\iso W_n$.
\item $\homp(\blank, W)$ takes colimits of simplicial sets to 
limits in $\cM \ovcat X$.  In other words, if $K = \colim_i K_i$, then
$\homp(K, W) \iso \lim^X_i \homp(K_i, W)$.
\end{enumerate}
\end{lemma}

\begin{defn}
\label{de:matching}
The object $\homp(\bd{n},W)$ is the \mdfn{$n$th augmented matching
space $M_n W$}.  The induced map $\homp(\del{n},W)\ra
\homp(\bd{n},W)$, which we may now write as $W_n \ra M_n W$, is the
\mdfn{$n$th matching map} for $W$.  Note that $W_0 \ra M_0 W$ is just
the augmentation since $\bd{0}=\emptyset$.
\end{defn}

We have chosen to work with these augmented constructions only because
they seem to make for the most compact and intuitive proofs.  Note
that the augmented matching objects and maps are the ones that arise
when considering Reedy model structures of simplicial objects in $(\cM
\ovcat X)$ \cite[15.2.2]{H}.  For $n\geq 2$, $M_n W$ is isomorphic to
$\tilde{M}_nU=\hom(\bd{n}, W)$ because $\bd{n}$ is connected.  The
following lemma is a reformulation of (\ref{eq:matchingmap}).

\begin{lemma}
\label{le:hc}
An augmented simplicial presheaf $U\ra X$ is a hypercover iff each
$U_n$ is a coproduct of repesentables and the maps $U_n \ra M_n U$
are all generalized covers.
\end{lemma}

\begin{defn}
\label{de:bndedhc}
A hypercover $U \ra X$ is \dfn{bounded} if there exists an $n\geq 0$
such that the maps $U_k \ra M_k U$ are isomorphisms for all $k> n$.
The smallest such $n$ for which this is true is called the
\dfn{height} of the hypercover, and denoted $\hgt U$.
\end{defn}

We have already remarked that the hypercovers of height $0$ are
precisely the \CCech complexes.  If one thinks of the $n$th level of a
hypercover as refining the $(n+1)$-fold `intersections' of the objects
in previous levels, then a bounded hypercover is one where the
refinement process stops at some point.  The following lemma is a
minor ingredient in the discussion of coskeleta in Section~\ref{se:cosk}
below, but the ideas from the proof reappear several times throughout
the paper.

\begin{lemma}
\label{le:fhc}
If $U\ra X$ is a bounded hypercover of height at most $n$, then the
induced maps $\homp(\del{k},U) \ra \homp(\sk_n\del{k},U)$ are
isomorphisms for all $k$.
\end{lemma}

\begin{proof}
When $k \leq n$, the result is easy because $\Delta^{k}$ equals $\sk_n
\del{k}$.  In general, $\del{k}$ is obtained from $\sk_n\del{k}$
by gluing on finitely many simplices of dimension at least $n+1$.
It suffices to show that $\homp(L,U) \ra \homp(K,U)$ is an isomorphism 
if $L$ is obtained from $K$ by attaching a simplex of dimension $i$,
where $i > n$.
Using Lemma~\ref{le:hompo} we obtain a pullback square
\[ \xymatrix{
 \homp(L,U) \ar[r]\ar[d] 
       & \homp(K,U) \ar[d]\\
 \homp(\del{i},U) \ar[r] & \homp(\bd{i},U).}
\]
The bottom map is the matching map $U_{i} \ra M_{i}U$, 
which is an isomorphism since $i > n$.
Hence the top map is also an isomorphism.
\end{proof}

\subsection{Skeleta, coskeleta, and split objects}\mbox{}\par
\label{se:cosk}

We continue to assume that $\cM$ is complete and cocomplete.  Let
$s\cM_{\leq n}$ and $s_+\cM_{\leq n}$ denote the categories of
$n$-truncated simplicial objects and augmented $n$-truncated
simplicial objects over $\cM$.  There is an obvious forgetful functor
$s_+\cM \ra s_+\cM_{\leq n}$ called $\sk_n$, and this has a right
adjoint called $\coskp_n$.  These are the \dfn{skeleta} and
\dfn{coskeleta} functors for augmented simplicial objects.  If $W$
belongs to $s_+\cM$, we abbreviate $\coskp_n \sk_n W$ to just
$\coskp_n W$.

The $k$th object of $\coskp_n U$ is
\[ [\coskp_n U]_k\iso\homp(\del{k},\coskp_n U)\iso\homp(\sk_n\del{k},U)
\]
(use Remark~\ref{re:adjoint} for the second isomorphism).  
In particular, the $(n+1)$st object of $\coskp_n U$ is what we have
been calling $M_{n+1}U$.  Observe also, using Lemma~\ref{le:fhc}, that
a hypercover $U$ has height at most $n$ if and only if $U\iso\coskp_n
U$.

\medskip

Now, the functor $\sk_n$ also has a left adjoint $\dgn_n\colon
s_+\cM_{\leq n} \ra s_+\cM$, called the \mdfn{$n$-degeneration}
functor.  The simplicial object $\dgn_n U$ is obtained from $U$ by
freely adding the images of the degeneracies in dimensions higher than
$n$ (and so, in particular, note that the augmentations are
irrelevant).  The object $[\dgn_n U]_{n+1}$ is called the
\mdfn{$(n+1)$st latching object} for $U$ and is denoted
\mdfn{$L_{n+1}U$}.  This latching object is the one that arises when
considering Reedy model structures of simplicial diagram categories
\cite[15.2.2]{H}.  Note that $U$, $\dgn_n U$, and $\cosk_n U$ all have
the same $n$-skeleton, so there are canonical maps
$\dgn_n U \ra U \ra \cosk_n U$; looking in level $n+1$ gives $L_{n+1}U
\ra U \ra M_{n+1}U$.

\begin{defn}
\label{de:split}
An object $W$ of $s\cM$ is said to be \dfn{split}, or to have
\dfn{free degeneracies}, if there exist subobjects $N_k\inc W_k$ such
that the canonical maps $N_k \amalg L_kW \ra W_k$ are isomorphisms for
all $k\geq 0$.  This is equivalent to requiring that the canonical map
\[ \coprod_{\sigma} N_\sigma \ra W_k
\]
is an isomorphism, where the variable $\sigma$ ranges over all
surjective maps in $\Delta$ of the form $[k]\ra [n]$, $N_\sigma$
denotes a copy of $N_n$, and the map $N_\sigma \ra W_k$ is the one
induced by $\sigma^*\colon W_n \ra W_k$ (see \cite[Def. 8.1]{AM}).
\end{defn}

The idea is that the objects $N_k$ represent the non-degenerate part
of $W$ in dimension $k$, and that the leftover degenerate part is as
free as possible.  The same definition as above can be applied to
augmented simplicial objects, and the result is that such an object is
split if and only if it is split when one forgets the augmentation.

We are particularly interested in \dfn{split hypercovers}.  If $U \map
X$ is a split hypercover then $L_k U$ is a summand of $U_k$, each $L_k
U$ is a coproduct of representables, and each representable summand of
$L_k U$ is the image under some degeneracy of a representable from
$U_{k-1}$ (but not uniquely).  It follows from \cite[Cor.~9.4]{D}
that split hypercovers are cofibrant in $U\cC_{\cL}$, which is why we
care about them.

\subsection{Computing matching objects}\mbox{}\par
\label{se:compmatch}

Suppose that $U\ra X$ is an augmented simplicial presheaf which in
each level is a coproduct of representables.  Note that (1) the
decomposition of $U_n$ into a coproduct of representables is unique up
to permutations of the summands, and (2) to give a map $\amalg_i A_i
\ra \amalg_j B_j$ between coproducts of representables corresponds to
giving, for each index $i$, an index $j(i)$ and a map $A_i \ra
B_{j(i)}$.  Because of these remarks, one can construct a simplicial
set $K$ by taking $K_n$ to be the set of representable summands of
$U_n$.  We'll refer to $K$ as the \dfn{indexing simplicial set} for
$U$.

Now suppose $a \colon L\ra K$ is a map of simplicial sets.  If
$\Delta^{op} L$ denotes the opposite category of simplices of $L$
\cite[Def. 15.1.16]{H}, there is an obvious diagram $\Delta^{op} L \ra
\sPre(\cC)\ovcat X$ which sends a $k$-simplex $\sigma$ to the
representable which is the summand of $U_k$ corresponding to
$a(\sigma)$.  We'll write $U(a)$ for the limit of this diagram.

The following observation is straightforward (use
Remark~\ref{re:adjoint}):

\begin{prop}
\label{pr:comp1}
There is an isomorphism of presheaves
\[\homp(L,U) \cong \coprod_{a\colon L\ra K}U(a).
\]  In particular, the matching object $M_nU$ is isomorphic to
$\coprod_{a\colon \bd{n} \ra K} U(a)$.
\end{prop}

Note that $\Delta^{op}L$ is an infinite category.  If $L$ has the
property that every nondegenerate simplex has nondegenerate faces
(e.g. $L=\bd{n}$), then one can use a smaller version.
Let $\Delta^{op}_{nd}L$ be the subcategory whose objects are the
non-degenerate simplices, and where the maps correspond to face maps.
Under the above assumption on $L$, it is an easy exercise to check
that $\Delta^{op}_{nd}L \inc \Delta^{op}L$ is final (use the fact that
in any simplicial set a degenerate simplex is an iterated degeneracy
of a unique nondegenerate simplex).  Hence the limit $U(a)$ can be
computed over $\Delta^{op}_{nd}L$ in practice.  

%%%%%%%%%%%%%%%%%%%%%%%%%%%%%%%%%%%%%%%%%%%%%%%%%%%%%%%%%%%%%%%%%%%%%%

\section{Hypercovers and lifting problems}
\label{se:hclift}

In Proposition~\ref{pr:loclift} we saw how local weak equivalences
relate to solutions of homotopy-lifting problems---one gets lifts
after passing from a representable to the elements of a covering
sieve.  Typically, these liftings can't be made compatible on the
different pieces of the sieve.  In this section we find that one {\it
can\/} arrange for this kind of compatibility by using hypercovers.

\medskip

The following proposition is the key ingredient in the proof of our
main result, Theorem \ref{th:hcloc}.  Recall that, just as for
ordinary covering families, a refinement of a hypercover $U \map X$ is
another hypercover $V \map X$ that factors through $U$.

\begin{prop}
\label{pr:hclift}
Let $F\ra G$ be a local acyclic fibration and let $K \map L$ be a
cofibration of finite simplicial sets.  For any square
\begin{myequation} 
\label{di:hclift}
\xymatrix{ K \tens U \ar[r]\ar[d] & F \ar[d] \\
              L \tens U \ar[r] & G
}
\end{myequation}
in which $U \map X$ is a hypercover, 
there exists another hypercover $V \ra X$ refining $U$
and liftings as in the following diagram: 
\begin{myequation} 
\label{di:hclift2}
\xymatrix{ 
K \tens V \ar[r] \ar[d] & K \tens U \ar[r] & F \ar[d] \\
L \tens V \ar[r] \ar@{.>}[urr] & L \tens U \ar[r] & G.
}
\end{myequation}
\end{prop}

To summarize the basic idea of the proof, 
%of Proposition~\ref{pr:hclift}, 
let's assume that $K\ra L$ is $\emptyset
\ra *$ and that the Grothendieck topology comes with a specified basis
of covering families.  Starting with a map $U\ra G$, we know by the
local-lifting property (\ref{pr:loclift}) that there is a covering
family $\{V_a \ra U_0\}$ with liftings $s_a \colon V_a \ra F$.  In
general, $s_a \restr{V_{ab}}$ and $s_b\restr{V_{ab}}$ are not equal,
but the two liftings become homotopic
after projecting down to $G$.  We can lift this homotopy to $F$ by
passing to a suitable covering family of $V_{ab}$, again using the
fact that $F\ra G$ is a local weak equivalence.  Next we move on to
consider patching on the triple intersections.  Once again, we can
patch up to homotopy after refining the triple intersections by a
covering family.  In this way we build a hypercover $V$ over which a
lifting is defined.  The work in this section is just a precise way of
saying all this.

\medskip

The proof involves an inductively constructed hypercover, and the
following lemma is the core of the induction step:

\begin{lemma}
\label{le:0lift}
Let $F$ and $G$ be presheaves of sets, and let $F\ra G$ be a
generalized cover.  If $J$ is a presheaf of sets with a map $J\ra
G$, then there exists a generalized cover $Z\ra J$ such that $Z$ is a
coproduct of representables and such that the diagram
\[ \xymatrix {&&F\ar[d] \\  Z\ar@{.>}[urr]\ar[r] & J \ar[r] &G
}
\]
has a lifting.
\end{lemma}

\begin{proof}
For each map $f\colon X\ra J$ from a representable, choose a covering
sieve $R_f$ of $X$ so that the composites $U\ra X\ra J \ra G$ lift to
$F$ for every $U\ra X$ in the sieve.  Let $Z$ denote the coproduct
\[ Z= \coprod_{X \llra{f} J} \Biggl ( \,
               \coprod_{\overset{U\ra X}{\text{in}\, R_f}} U
              \Biggr ).
\]
The obvious map $Z\ra J$ is a generalized cover, and the
composite $Z\ra J\ra G$ lifts to $F$.
\end{proof}

\begin{prop}
\label{pr:hcind}
Let $F\ra G$ be a local acyclic fibration, 
and let $U \map G$ be a map where $U \map X$ is a hypercover.
Let $n \geq 0$, and suppose that there is
an $n$-truncated hypercover $V \ra X$ refining $\sk_n U$ and a map
$V\ra F$ such
that the diagram
\begin{myequation}
\label{eq:extend}
\xymatrix{ && F \ar[d] \\
           V \ar[r]\ar[urr] & U \ar[r] & G
}
\end{myequation}
commutes.  Then there is an $(n+1)$-truncated hypercover $W \ra X$
refining $U$ and a map $W \ra F$ making the corresponding diagram
commute, and such that on $n$-skeleta the diagram is equal to
(\ref{eq:extend}).
\end{prop}

\begin{proof}
The core of the proof is just an Artin-Mazur argument
\cite[Ch. 8]{AM}.  First form the pullback $F'=U\times_G F$.  The map
$F'\ra U$ is still a local acyclic fibration, and we need only produce
an $(n+1)$-truncated hypercover $W$ and a lifting into $F'$.  In other
words, we can reduce to the case where $U=G$ (and $F=F'$).  Note that
in this case $G$ is locally fibrant---the representable $X$ is locally
fibrant for trivial reasons, and $U\ra X$ is a local fibration.
Moreover, since $F\ra G$ is a local fibration, $F$ is locally fibrant
as well.

Now $F^{\del{n+1}}\ra F^{\bd{n+1}}$ is a local fibration by
Proposition~\ref{pr:fsm7}, so the map in the $0$th level is a
generalized cover by (\ref{eq:matchingmap}).  When $n>0$ this map is
precisely $F_{n+1} \map M_{n+1}F$ (the $n=0$ case being only
slightly different).
Our initial diagram gives a map
$V^{\bd{n+1}} \ra F^{\bd{n+1}}$, and the $0$th level has the form
$M_{n+1}V \map M_{n+1} F$.  So Lemma~\ref{le:0lift} says that there
is a generalized cover $Z\ra M_{n+1}V$, where $Z$ is a coproduct of
representables, such that the composite $Z\ra M_{n+1}F$ lifts through
$F_{n+1}$.  We take $W$ to be the $(n+1)$-truncated hypercover
with $\sk_n W=\sk_n V$ and $W_{n+1}=Z\amalg L_{n+1} V$.
\end{proof}

\begin{proof}[Proof of Proposition~\ref{pr:hclift}]
Given a square as in the statement of the proposition, it may be
interpreted as a map
\[ U \ra F^{K}\times_{G^{K}} G^{L}. \]
We are trying to produce a hypercover $V\ra X$ refining $U \map X$ 
and a lifting
\[
\xymatrix{ && F^{L} \ar[d] \\
           V \ar[r]\ar@{.>}[urr] & U \ar[r] 
                  & F^{K}\times_{G^{K}} G^{L}.
}
\]
The vertical map is a local acyclic fibration by Proposition~\ref{pr:fsm7},
so the hypercover can be produced inductively using
Proposition~\ref{pr:hcind}.
\end{proof}

If we have a map $F\ra G$ which is a local weak equivalence but not
necessarily a fibration, we can say the following:

\begin{prop}
\label{pr:hclift2}
Let $F \ra G$ be a local weak equivalence between locally fibrant
simplicial presheaves.  Then given any diagram as in
(\ref{di:hclift}), there exists a hypercover $V \map X$ refining $U
\map X$ and relative-homotopy-liftings in the diagram
(\ref{di:hclift2}).
\end{prop}

Recall that relative-homotopy-liftings were defined in
Proposition~\ref{pr:loclift}, and discussed extensively in \cite{DI2}.

\begin{proof}
Given a diagram as in (\ref{di:hclift}), we need to produce a
hypercover $V \map X$ refining $U \map X$ together with liftings in
the diagram
\[
 \xymatrixcolsep{1.5pc}\xymatrix{
  &K \tens V \ar[r]\ar[d]\ar[dl] & F \ar[d] \\
  L\tens V \ar[d]_{i_0} \ar@{.>}[urr]  
             & L \tens V \ar[r]\ar[dl]_{i_1} & G  \\
  RH \tens V.  \ar@{.>}[urr]
}
\]
Here $RH$ denotes the pushout of $L\times \del{1} \lla K \times\del{1}
\llra{\pi} K$, and the maps $i_0$ and $i_1$ are the obvious inclusions
$L\inc RH$.

Consider the square
\[ \xymatrix{ 
    F^{RH} \ar[d]_{i_1^*}\ar[r] 
           & F^{L}\times_{G^{L}} G^{RH} \ar[d]\\
    F^{L} \ar[r] 
           & F^{K}\times_{G^{K}} G^{L}.
}
\]
By \cite[Cor. 7.5]{DI2}, the fact that $F\ra G$ is a local weak
equivalence between locally fibrant objects implies that the
horizontal maps are also local weak equivalences.  By the same result,
the fact that $i_1\colon L \ra RH$ is a weak equivalence of simplicial
sets implies that the left vertical map is a local weak equivalence.
So we conclude that the same is true of the right vertical map.  Even
more, the right vertical map is a local fibration by
Lemma~\ref{le:cornerfib} below.

Our initial data from (\ref{di:hclift}) was a map $U \ra
F^{K}\times_{G^{K}} G^{L}$, so by Proposition~\ref{pr:hclift} (for
$n=0$ and $K\ra L$ equal to $\emptyset \ra *$) there is a hypercover
$V\ra X$ refining $U \map X$ for which the composite lifts through
$F^{L}\times_{G^{L}} G^{RH}$.  This provides the necessary
relative-homotopy-lifting.
\end{proof}

\begin{lemma}
\label{le:cornerfib}
Let $F\ra G$ be a map between locally fibrant simplicial presheaves.
Assume we have a square of finite simplicial sets
\[\xymatrix{ K\ar[d]\ar[r] & M \ar[d] \\
                L\ar[r]  & N}
\]
such that both $K\ra M$ and $M\amalg_K L \ra N$ are
cofibrations.  Then the induced map
\[ F^M\times_{G^M} G^N \ra F^K \times_{G^K} G^L \]
is a local fibration.
\end{lemma}

\begin{proof}
The hypotheses imply that both $F^M\ra F^K$ and $G^N \ra G^{(M\amalg_K
L)}=G^M\times_{G^K} G^L$ are local fibrations, using
\cite[Cor. 1.5]{J2}.  Now we observe that there are pullback squares
\[ \xymatrix{
F^M \times_{G^M} G^N \ar[r]\ar[d] & G^N \ar[d] 
           && F^M\times_{G^K} G^L\ar[r]\ar[d]& F^M \ar[d]\\
F^M \times_{G^K} G^L \ar[r] & G^M\times_{G^K} G^L 
           && F^K\times_{G^K} G^L \ar[r] &F^K,
}
\]
and the pullback of a local fibration is again a local fibration.
Finally, the map we want is just the composite $F^M\times_{G^M} G^N
\ra F^M\times_{G^K}G^L \ra F^K\times_{G^K} G^L$.
\end{proof}

%%%%%%%%%%%%%%%%%%%%%%%%%%%%%%%%%%%%%%%%%%%%%%%%%%%%%%%%%%%%%%%%%%%%

\section{Hypercovers and localizations}
\label{se:hcloc}

In this section we prove the main theorem, that Jardine's model
category can be obtained by localizing the injective structure
$\sPre(\cC)$ at the hypercovers.  This lets us identify
the fibrant objects in the model structure.  Similar results are
proven for the projective version $U\cC_{\cL}$.

\medskip

We start with a definition:

\begin{defn}
A collection of hypercovers $S$ is called \dfn{dense} if every
hypercover $U\ra X$ in $\sPre(\cC)$ can be refined by a hypercover
$V\ra X$ which belongs to $S$.  
%The collection is called \dfn{split}
%if every hypercover in $S$ can be refined by a split hypercover which
%also belongs to $S$.
\end{defn}

For instance, Theorem~\ref{th:bhcwc} shows that when $\cC$ is a
Verdier site the collection of basal hypercovers is both split and
dense.  The following is our main goal.

\begin{thm}
\label{th:hcloc}
Let $S$ be a collection of hypercovers which contains a set that is
dense (e.g., the collection of all hypercovers).  
Then the localization $\sPre(\cC)/S$ exists and coincides with
Jardine's model structure $\sPre(\cC)_{\cL}$.  
Similarly, the localization $U\cC/S$ exists and coincides
with $U\cC_{\cL}$.
\end{thm}

Our notation is that if $\cM$ is a model category and $S$ is a
collection of maps, then $\cM/S$ denotes the left Bousfield
localization of $\cM$ at $S$ (if it exists)---see \cite{D} for a
summary treatment or \cite{H} for complete details.  The fibrations
and weak equivalences in $\cM/S$ are called $S$-fibrations and
$S$-equivalences, while the cofibrations are the same as those in
$\cM$.

The hypothesis of the theorem is a little stronger than just assuming
that $S$ is dense, because $S$ itself may not be a set. 
For the same reason, 
the existence of the localization is not automatic.
One of the things we will do is apply this theorem in the case where
$S$ is the collection of {\it all\/} hypercovers, and this is not a
set: in our 
definition of hypercover one can have
arbitrarily large coproducts of representables appearing.  
So we'll need to verify that $S$ contains a dense {\it set\/}, and
this can be done by
making use of the fact that
our site is small. We choose a suitably large regular cardinal, and
then we only consider hypercovers in which the number of summands in
each level is bounded by our cardinal.  For now we can ignore this
point, but see Section~\ref{se:card}.

To prove Theorem \ref{th:hcloc}, we need a general criterion for
checking whether two localizations are identical:

\begin{lemma}
\label{le:locfact} 
Let $\cM$ be a model category, and let $S \subseteq T$ be
two classes of maps for which the localizations $\cM/S$ and $\cM/T$ exist.
For the two localizations to be the same it suffices to check the
following: if an $S$-fibration $X\fib Y$ between $S$-fibrant objects
is a $T$-equivalence, then it is an $S$-equivalence.
\end{lemma}

\begin{proof}
We must show the hypothesis implies that every $T$-equivalence $A\ra
B$ is an $S$-equivalence.  Let $L$ denote a fibrant replacement
functor in $\cM/S$, and consider the square
\[ \xymatrix{ A \ar[r]\ar[d]_{\sim}^S &B\ar[d]_{\sim}^S \\
              LA \ar[r] & LB.  }
\] 
Since $S\subseteq T$ the two vertical maps are $T$-equivalences, and
the top map is a $T$-equivalence by assumption---so the bottom map is
one as well.  Now factor the bottom map in $\cM/S$ as an $S$-acyclic
cofibration followed by an $S$-fibration:
\[
\xymatrix{ LA \ar@{ >->}[r]_{\sim}^S & X \ar@{->>}[r] &LB.
}
\] 
Note that $X$ is $S$-fibrant, because $LB$ is.  Also, since both
the first map and the composite are $T$-equivalences, so is the second
map.

Therefore the map $X\ra LB$ is a $T$-equivalence and an $S$-fibration,
and the domain and codomain are
$S$-fibrant.  Our hypothesis then says that $X\ra LB$ is an
$S$-equivalence.  Applying the two-out-of-three property (twice)
shows that $A\ra B$ is an $S$-equivalence.
\end{proof}

For the moment let $S$ be a {\it set\/} of hypercovers that is dense.
Because $S$ is a set, we know that the model structure $\sPre(\cC)/S$
exists (by \cite[Thm. 4.1.1]{H}, using that $\sPre(\cC)$ is left
proper and cellular).  The fibrant objects in $\sPre(\cC)/S$ (called
$S$-fibrant objects) are the injective-fibrant objects which satisfy
descent for all hypercovers in $S$.  Since every hypercover is a local
weak equivalence by definition, $\sPre(\cC)_{\cL}$ is a localization
of $\sPre(\cC)/S$.  To show that the two structures coincide, we now
check the criterion from the above lemma:

\begin{lemma}
\label{le:heart}
Let $F$ and $G$ be $S$-fibrant objects, and let $f:F\ra G$ be an
$S$-fibration that is also a local weak equivalence.  If
$X$ is a representable, then every square
\[
\xymatrix{
\bd{n} \tens X \ar[r] \ar[d] & F \ar[d] \\
\del{n} \tens X \ar[r] & G     }
\]
has a lifting.
In particular, $f$ is actually an objectwise acyclic fibration
and therefore an $S$-equivalence.
\end{lemma}

\begin{proof}
The second claim follows from the first by adjointness and because
acyclic fibrations of simplicial sets are detected by the 
right lifting property with respect to the maps $\bd{n} \map \del{n}$.

Now we prove the first claim.
First, 
$f$ is an objectwise fibration since 
every $S$-fibration is an injective-fibration and also a projective-fibration.
This implies that $f$ is also a local fibration.
Because $f$ is both a local fibration and a local weak equivalence,
Proposition~\ref{pr:hclift} guarantees us a hypercover $U \ra X$
such that the diagram
\[ 
\xymatrix{ \bd{n} \tens U \ar[r] \ar[d] & \bd{n} \tens X \ar[r] & F
\ar[d] \\ \del{n} \tens U \ar[r] \ar@{.>}[urr] & \del{n} \tens X
\ar[r] & G }
\] 
has a lifting.  In applying Proposition \ref{pr:hclift}, we have used
that $X$ is (trivially) a hypercover of itself.  Since $S$ is dense,
we may refine $U$ and assume that $U \map X$ belongs to $S$.  We now
write down the following diagram of simplicial mapping spaces:
\[
\xymatrix{ \Map(X,F^{\del{n}}) \ar[r]^\sim \ar@{->>}[d] & 
               \Map(U,F^{\del{n}}) \ar@{->>}[d] \\
           \Map(X,G^{\del{n}} \times_{G^{\bd{n}}} F^{\bd{n}}) \ar[r]^\sim & 
  \Map(U,G^{\del{n}} \times_{G^{\bd{n}}} F^{\bd{n}}). 
}
\]
All the model categories we have been considering are simplicial model
categories, and this implies that $F^{\del{n}} \map G^{\del{n}}
\times_{G^{\bd{n}}} F^{\bd{n}}$ is an $S$-fibration between
$S$-fibrant objects.  
%It is a local weak equivalence by
%Proposition \ref{pr:fsm7}.  
Therefore, the vertical maps above are
fibrations of simplicial sets because both $X$ and $U$ are
$S$-cofibrant.  Likewise, the horizontal maps are weak
equivalences because $U \ra X$ is an $S$-equivalence between
$S$-cofibrant objects and both $F^{\del{n}}$ and $G^{\del{n}}
\times_{G^{\bd{n}}} F^{\bd{n}}$ are $S$-fibrant.

We are given a $0$-simplex $x$ in the lower left corner in the above
diagram, and we want to find a
lift in the upper left corner.  We have already shown that the image of
$x$ in the lower right corner lifts to the upper right corner.
Since the horizontal maps are weak equivalences,
there is another $0$-simplex $y$ belonging to the connected component
of $x$ such that $y$ has a lift in the upper left corner.
But fibrations of simplicial sets are
surjective onto the components in their images, so $x$ also has a lift.
\end{proof}

\begin{proof}[Proof of Theorem~\ref{th:hcloc}]
We first consider the claim for $\sPre(\cC)_{\cL}$.  For the case when
our collection of hypercovers $S$ is itself a {\it set\/}, we have
already done all the work.  Since hypercovers are local weak
equivalences we know $S\subseteq \cL$, and so we are in the situation
of Lemma~\ref{le:locfact}.  The necessary condition was verified in
Lemma~\ref{le:heart}.

In the general case, let $S'$ be a dense set of hypercovers contained
in $S$.  As shown in the previous paragraph, $\sPre(\cC)/S'$ is equal
to $\sPre(\cC)_{\cL}$.  So every local weak equivalence is a weak
equivalence in $\sPre(\cC)/S'$, and in particular every hypercover in
$S$ is an $S'$-equivalence.  This shows that $\sPre(\cC)/S$ exists and
is equal to $\sPre(\cC)/S'$.

The argument for $U\cC/S$ is basically the same.  Assume first that
$S$ is a set of hypercovers which is dense.  One reproves the analog
of Lemma~\ref{le:heart} for the projective model structure; the only
difference in the proof is that one replaces $U$ by a cofibrant object
before dealing with simplicial mapping spaces.  The rest of the
argument is exactly the same, as is the generalization to the case
where $S$ need not be a set.
\end{proof}

\subsection{Cardinality considerations}
\label{se:card}

Early in this section we mentioned that the collection of all
hypercovers is not a set, but contains a subset that is dense.  We
will now give the proof.  Recall from Section \ref{se:compmatch} that to any
hypercover $U\ra X$ one can attach an indexing simplicial set $K$,
where $K_n$ is the set of representable summands of $U_n$.  The
\dfn{size} of the hypercover is the cardinality of $\coprod_n K_n$,
i.e., the number of representable summands that appear in $U$.  The
main point is that in the arguments from Proposition~\ref{pr:hcind}
and Lemma~\ref{le:0lift}, one can control the size of the constructed
hypercover.

\begin{prop}
\label{pr:denseset}
The class of all hypercovers has a subset which is dense.
\end{prop}

\begin{proof}
Choose a regular cardinal $\lambda$ sufficiently large compared to the
cardinality of the set of morphisms in $\cC$, and let $S$ denote the set
of all hypercovers of size less than $\lambda$.  We will show that any
hypercover $U\ra X$ can be refined by one in $S$.

Since $U_0\ra X$ is a generalized cover, there is a covering sieve $R$
of $X$ such that every $W\ra X$ in $R$ lifts through $U_0$.  Let
$V_0=\coprod_{W\ra X} W$, where the coproduct ranges over all maps
$W\ra X$ in $R$.  The number of summands in $V_0$ is clearly bounded
by $\lambda$.

Now assume by induction that we have constructed an $n$-truncated
hypercover $V\ra X$ which refines $U$, and such that the number of
summands in $V$ is less than $\lambda$.  To extend $V$ we use the
argument from Proposition~\ref{pr:hcind}, where we must show that
$Z$ does not have too many representable summands.  
Inspecting the construction of $Z$ given in Lemma~\ref{le:0lift},
it suffices to show that there aren't too many maps from a representable
into $M_{n+1} V$.
This can be deduced from Proposition \ref{pr:comp1}.
\end{proof}

%%%%%%%%%%%%%%%%%%%%%%%%%%%%%%%%%%%%%%%%%%%%%%%%%%%%%%%%%%%%%%%%%%%%

\section{Fibrations and descent conditions}
\label{se:fibs}

This section identifies the fibrations and fibrant objects in
$\sPre(\cC)_{\cL}$ and $U\cC_{\cL}$.  We start with the fibrant
objects, where the result follows from Theorem~\ref{th:hcloc}:

\begin{cor}
\label{co:hcloc}
Let $S$ be a collection of hypercovers which contains a set that is
dense.  A simplicial presheaf $F$ is fibrant in $\sPre(\cC)_{\cL}$
(resp. in $U\cC_{\cL}$) if and only if $F$ is injective-fibrant
(resp. objectwise-fibrant) and satisfies descent for all hypercovers
in $S$.  
\end{cor}

Note that an immediate consequence is that a simplicial
presheaf $F$ satisfies descent for all hypercovers if and only if it
satisfies descent for all elements of $S$.

\begin{proof}
First observe that the fibrant objects in $\sPre(\cC)/S$ are the
injective-fibrant objects $F$ such that $\Map(X, F) \map \Map(U, F)$ is a
weak equivalence for every $U\ra X$ in $S$ \cite[Thm.~4.1.1(2)]{H}.
(Since everything is cofibrant in $\sPre(\cC)$, one doesn't have to
take cofibrant replacements for $U$ and $X$.)
Lemma~\ref{le:descent}(i) says the latter condition is the same as $F$
satisfying descent for $U$.  By Theorem~\ref{th:hcloc} the model
structure $\sPre(\cC)_{\cL}$ is the same as $\sPre(\cC)/S$, so this
proves one case.

The proof of the second case is the same, but any hypercover $U$ must be
replaced by a cofibrant object before it appears in a mapping space, and
Lemma~\ref{le:descent}(ii) is used instead of Lemma~\ref{le:descent}(i).
\end{proof}

Using the above corollary, we can actually identify the fibrations in
both $\sPre(\cC)_{\cL}$ and $U\cC_{\cL}$ in terms of a relative
descent condition.  The idea for this proof is due to Blander
\cite[Prop. 4.1]{Bl}.  However, the interpretation in terms of descent
conditions requires extra care with various homotopy-limit
constructions.  We start by recalling the relevant definitions.

Let $P$ be the cosimplicial simplicial set $B(\Delta\ovcat \blank)$
described in \cite[14.7.7]{H}.  Then $\hom^\Delta(P,X)$ is
defined in \cite[18.3.6]{H} to be the homotopy limit of a cosimplicial
diagram $X$.  However, this construction only has good properties when
$X$ is objectwise-fibrant.  Therefore, we will define $\holim_\Delta X$
to be $\hom^\Delta(P, \hat{X})$, where $\hat{X}$ is an objectwise-fibrant
replacement for $X$.  If $X$ is not objectwise-fibrant, there is no
obvious guarantee that $\hom^\Delta(P,X)$ and $\holim_\Delta X$ are
weakly equivalent; we will need to be careful about this below.  

In the following definition, $F(U)$ is the cosimplicial simplicial set
which appears in Definition~\ref{de:hc-desc}.  If $Z$ is any
simplicial set, then $cZ$ denotes the corresponding {\it constant\/}
cosimplicial simplicial set.

\begin{defn}
\label{defn:relative-descent}
An objectwise fibration $F\ra G$
\dfn{satisfies descent} for the hypercover $U\ra X$ if the natural map 
$F(X)\ra \holim_\Delta [cG(X)\times_{G(U)} F(U)]$ is a weak equivalence.
\end{defn}

Note that this reduces to our previous definition when $G=*$.

\begin{remark}
By manipulating homotopy limits, one can see that the above definition
is equivalent to requiring that
$F(X)$ be the homotopy limit of the diagram
\[ \xymatrix{
&   \prod_a F(U^a_0) \ar@<0.5ex>[r]\ar@<-0.5ex>[r]\ar[d] &
    \prod_{a} F(U^a_{1}) \ar@<0.6ex>[r]\ar[r]\ar@<-0.6ex>[r]\ar[d] &  \cdots \\
G(X) \ar[r] &   \prod_a G(U^a_0) \ar@<0.5ex>[r]\ar@<-0.5ex>[r] &
    \prod_{a} G(U^a_{1}) \ar@<0.6ex>[r]\ar[r]\ar@<-0.6ex>[r] &  \cdots \\
}
\]
or alternatively of the diagram
\[
G(X) \map \holim_\Delta G(U) \leftarrow \holim_\Delta F(U).
\]
We will not need either of these criteria, however.
\end{remark}

Our goal is to show that the fibrations in $\sPre(\cC)_{\cL}$ or in
$U\cC_{\cL}$ can be characterized using the above descent condition.  The
proof is more complex than one might imagine, and proceeds in a few
steps.

\begin{lemma}\mbox{}\par
\label{le:holim-criterion}
\begin{enumerate}[(i)]
\item
Let $Z$ be a simplicial set, let $W^*$ be a cosimplicial simplicial set, 
and suppose there is an
objectwise fibration $W^* \ra cZ$.  Then $\hom_\Delta(P,W)$ is weakly
equivalent to the homotopy limit of  $W^*$.
\item An objectwise fibration of simplicial presheaves $F\ra G$
satisfies descent for $U \ra X$ if and only if the natural map $F(X) \ra
\hom^\Delta(P, cG(X)\times_{G(U)} F(U))$ is a weak equivalence.
\end{enumerate}
\end{lemma}

\begin{proof}
For (i), pick a fibrant
replacement $Z \ra \hat{Z}$, and factor the composite $W \ra cZ \ra
c\hat{Z}$ into an objectwise acyclic cofibration $W \trcof \hat{W}$
followed by an objectwise fibration $\hat{W} \fib c\hat{Z}$.  Let
$E=(cZ)\times_{c\hat{Z}} \hat{W}$, and consider the two maps
$\hom_{\Delta}(P,W) \ra \hom_{\Delta}(P,E) \ra
\hom_{\Delta}(P,\hat{W})$.  The last object is the homotopy limit of
$W$ by definition, so we will show that both maps are equivalences.

First consider the diagram
\[ \xymatrix{
Z \ar[r]^-\sim\ar[d]_\sim &\hom_{\Delta}(P,cZ) \ar[d] 
      &\hom_{\Delta}(P,E)\ar[d]\ar[l] \\
\hat{Z} \ar[r]^-\sim &\hom_{\Delta}(P,c{\hat Z}) 
      &\hom_{\Delta}(P,\hat{W})\ar@{->>}[l]. 
}
\]
Note that $P_0 = B(\Delta \ovcat [0])$ is actually simplicially
contractible, since $[0]$ is the terminal object of $\Delta$.  Because
$\hom^{\Delta}(P,cZ) =\Map(P_0,Z)$ it follows that the canonical map
$Z=\hom^{\Delta}(*,cZ) \ra \hom^{\Delta}(P,cZ)$ is a simplicial
homotopy equivalence (and similarly for $\hat{Z}$).  So the left
horizontal maps in the diagram are simplicial homotopy equivalences.

The right square is a pullback square.  It follows that the right-most
vertical map is the pullback of a weak equivalence along a fibration,
hence a weak equivalence (by right properness of $\sSet$).

Next we consider the diagram 
\[ \xymatrix{W \ar[rr]^\sim \ar@{->>}[rd] && E\ar@{->>}[ld] \\
 & cZ.
}
\]
It is a general fact that a right Quillen functor preserves weak
equivalences between fibrations over a given object (this is Ken
Brown's Lemma \cite[7.7.2]{H} in the overcategory).  Since
$\hom_{\Delta}(P,\blank)$ is a right Quillen functor, we therefore
have $\hom_{\Delta}(P,W) \we \hom_{\Delta}(P,E)$.  This finishes the
proof of (i).

Part (ii) is an application of (i), noting that $cG(X)\times_{G(U)}
F(U)\ra cG(X)$ is an objectwise fibration.  
\end{proof}

Our next task is to reinterpret descent in yet another way, as a
lifting condition.  In $U\cC$, consider the maps $\hocolim_n U_n \ra
X$ and factor them as $\hocolim_n U_n \cof B(U) \trfib X$, where the
first map is a projective cofibration and the second is an objectwise
acyclic fibration.  From now on we'll denote the object $\hocolim_n
U_n$ by $A(U)$, for short.  Note that since each $U_n$ is
projective-cofibrant, so is $\hocolim U_n$ by \cite[18.4.2]{H}---so
both $A(U)$ and $B(U)$ are projective-cofibrant.  Also, $X$ is both
cofibrant and fibrant in $U\cC$, therefore $B(U)$ is cofibrant-fibrant
and $B(U)\ra X$ is a simplicial homotopy equivalence.

Fix a collection of hypercovers $S$ which contains a dense set.  
Let $J_{\cC}$ be the collection of all maps $\Lambda^{n,k}\tens Z \ra
\Delta^n \tens Z$ (for all $Z \in \cC$) and also of all the maps
\[ [A(U)\tens \Delta^n] \amalg_{A(U)\tens \bd{n}} [B(U)\tens
\bd{n}] \ra B(U)\tens\del{n}
\]
where $U\ra X$ ranges over the hypercovers in $S$. 

\begin{lemma}
\label{le:descent=lift}
A map $F\ra G$ is an objectwise fibration satisfying descent with
respect to all hypercovers in $S$ if and only if it has the
right-lifting-property with respect to the maps in $J_{\cC}$.
\end{lemma}

Note that, as a consequence, objectwise fibrations satisfying descent
are closed under pullbacks.

\begin{proof}
A map $F\ra G$ is an objectwise fibration if and only if it has the
right-lifting-property with respect to the maps $\Lambda^{n,k}\tens Z
\ra \Delta^n \tens Z$.  If $F\ra G$ is an objectwise fibration then
the map 
\[
\Map(B(U),F) \ra \Map(A(U),F)\times_{\Map(A(U),G)} \Map(B(U),G)
\] 
is a fibration because $U\cC$ is a simplicial model
category, and it's a weak equivalence if and only if $F\ra G$ has the
right-lifting-property with respect to the maps 
\[
[A(U)\tens \Delta^n] \amalg_{A(U)\tens \bd{n}} [B(U)\tens \bd{n}] \ra 
    B(U)\tens\del{n},
\]
for all $n$.  Because $B(U) \ra X$ is a simplicial homotopy
equivalence, it follows that both $\Map(X,F) \ra \Map(B(U),F)$ and
$\Map(X,G)\ra \Map(B(U),G)$ are simplicial homotopy equivalences.
From this one sees that 
\[
\Map(B(U),F) \ra \Map(A(U),F)\times_{\Map(A(U),G)} \Map(B(U),G)
\]
is a weak equivalence
if and only if 
\[
\Map(X,F) \ra \Map(A(U),F)\times_{\Map(A(U),G)} \Map(X,G)
\]
is one.

We have $A(U)=\hocolim_{\Delta^{op}} U =U
\tens_{\Delta} P$, which gives the canonical identification
of $\Map(A(U),F)$ with $\hom^{\Delta}(P,F(U))$ \cite[18.1.10]{H}.
A diagram chase now shows that an objectwise fibration $F\ra G$
has the lifting property in question if and only if the natural maps 
\[ \hom^\Delta(P,cF(X)) \ra \hom^{\Delta}(P,F(U))
\times_{\hom^{\Delta}(P,G(U))} \hom^{\Delta}(P,cG(X)) 
\]
are weak equivalences.  The object on the left is simplicially
homotopy equivalent to $F(X)$, as in the proof of Lemma
\ref{le:holim-criterion}.  Using that $\hom^\Delta(P,\blank)$ is a
right adjoint, the object on the right can be identified with
$\hom^\Delta(P,cG(X)\times_{G(U)} F(U))$.  Now
Lemma~\ref{le:holim-criterion}(ii) tells us that the above map is a
weak equivalence if and only if $F \map G$ satisfies descent for $U
\map X$.
\end{proof}

\begin{lemma}
\label{le:j-inj}
Let $\cM$ be a model category, and let $J$ be a set of acyclic
cofibrations which permits the small object argument.  Suppose 
every map that is both a $J$-injective and a weak equivalence is also
a fibration.  Then the $J$-injectives are precisely the fibrations.
\end{lemma}

\begin{proof}
The small object argument shows that any acyclic cofibration $A\ra B$
may be factored as $A \ra X \ra B$ where the first map is a relative
$J$-cell complex (therefore an acyclic cofibration) and the second is
a $J$-injective.  The two-out-of-three property says that $X\ra B$ is
a weak equivalence, and then our assumption implies it is a fibration.
The retract argument then shows that $A\ra B$ is a retract of $A\ra
X$.  In other words, every acyclic cofibration is a retract of a
relative $J$-cell complex.  From this it follows that $J$-injectives
have the right-lifting-property with respect to all acyclic
cofibrations, hence $J$-injectives are fibrations.
\end{proof}

\begin{prop}
\label{pr:trfibs}
The set $J_{\cC}$ is a set of generating acyclic cofibrations in
$U\cC_{\cL}$.
\end{prop}

\begin{proof}
The maps in $J_{\cC}$ are acyclic cofibrations in
$U\cC_{\cL}$ because $U\cC_{\cL}$ is a simplicial model category and
because the maps $A(U) \map B(U)$ are acyclic cofibrations in
$U\cC_{\cL}$.
It suffices to show that if $F\ra G$ is a
$J_{\cC}$-injective which is also a local weak equivalence, then $F\ra
G$ is an objectwise acyclic fibration.  This will imply that it is
also an acyclic fibration in $U\cC_{\cL}$, and in particular a
fibration in $U\cC_{\cL}$.  Then we apply Lemma~\ref{le:j-inj}.

We need to show that for every $X\in \cC$ and every point $x\in G(X)$,
the fiber of $F(X) \ra G(X)$ over $x$ is contractible.  Let us replace
$F$ and $G$ by their restrictions to the site $\cC\ovcat X$.  The map
of restricted presheaves is still a local weak equivalence, and still
satisfies descent with respect to a dense set of hypercovers (because
every hypercover in $\cC\ovcat X$ is essentially a hypercover in
$\cC$).  By Lemma \ref{le:descent=lift} $F\ra G$ has the
right-lifting-property with respect to the corresponding set
$J_{(\cC\ovcat X)}$ in $U(\cC\ovcat X)$.

Our point $x\in G(X)$ now corresponds to a map $* \ra G$ in
$U(\cC\ovcat X)$.  Consider the pullback square
\[ \xymatrix{
         H \ar[r]\ar[d] & F\ar[d] \\
          {*} \ar[r]  & G.
}
\]
The map $H \ra *$ still has the right-lifting-property with respect to
$J_{(\cC\ovcat X)}$, and so it is an objectwise fibration satisfying
descent by Lemma \ref{le:descent=lift}.  
Therefore $H$ is fibrant in $U(\cC\ovcat X)_{\cL}$ by
Corollary~\ref{co:hcloc}.  Moreover, since $F\ra G$ is a local
acyclic fibration, so is $H \ra *$.  Thus $H \ra *$ is an acyclic
fibration in $U(\cC\ovcat X)_{\cL}$, hence an objectwise acyclic
fibration.  This implies that $H(*)$ is contractible.  But $H(*)$
is another name for the fiber of our original map $F(X) \ra G(X)$, so
we are done.
\end{proof}

\begin{thm}
\label{th:fibrations}
Let $S$ be a collection of hypercovers which contains a set that is
dense.
A map of simplicial presheaves $F\ra G$ is a fibration in 
$\sPre(\cC)_{\cL}$ (resp. in $U\cC_{\cL}$) if and only
if it is an injective fibration 
(resp. an objectwise fibration)
and satisfies descent for all hypercovers in $S$.  
\end{thm}

\begin{proof}
The statement for $U\cC_{\cL}$ follows from
Proposition~\ref{pr:trfibs} and Lemma~\ref{le:descent=lift}.  For
$sPre(\cC)_{\cL}$ one repeats all the above arguments, but in the
definition of $J_\cC$ the maps $\Lambda^{n,k}\tens Z \ra
\Delta^{n}\tens Z$ are replaced with a set of generating acyclic
cofibrations for $\sPre(\cC)$.
\end{proof}

\subsection{A short example about fibrant replacement}
\label{se:shfcohom}

We end this section with a simple (and well-known) example
demonstrating the use of Corollary~\ref{co:hcloc}.  Let $A$ be a
presheaf of abelian groups on the site $\cC$, and let $I_*$ denote an
injective resolution of the sheafification
$\tilde{A}$ in the category of sheaves.  We will
explain how to use $I_*$ to construct a fibrant replacement for the
simplicial presheaf $K(A,n)$.

Let $\cI$ denote the chain complex of presheaves which has $I_k$ in
dimension $n-k$ when $k<n$, and has the presheaf of $n$-boundaries
$B_n$ in dimension $0$.  The Dold-Kan correspondence lets us identify
presheaves of (non-negatively graded) chain complexes with the abelian
group objects in $\sPre(\cC)$, and so $\cI$ can be regarded as a
simplicial presheaf.  Since right now we are only dealing with abelian
things, it's easier just to think about chain complexes, though.

The map $A\ra I_0$ induces a map $K(A,n) \ra \cI$ (and recall that as
a chain complex, $K(A,n)$ has $A$ in dimension $n$ and $0$ everywhere
else).  This map is a local weak equivalence because it induces
isomorphisms on homology group sheaves.  We claim that $\cI$ satisfies
descent for all hypercovers, and so is a fibrant object in
$U\cC_{\cL}$.  This shows that one can identify weak homotopy classes
of maps $\Ho(\Delta^k/\bd{k}\tens X,K(A,n))$ with $H_{k}(\cI(X))$,
which is just the sheaf cohomology group $H_{\cC}^{n-k}(X,\tilde{A})$.
The connection with sheaf cohomology is also explained in
\cite[Section 2]{J1}.

If $U\ra X$ is a hypercover of $X$, let $\Z[U]$ denote the chain
complex of presheaves obtained by applying the free abelian group
functor to the presheaves $U_n$.  It is known that after
sheafification $\Z[U]$ becomes a resolution of $\Z[X]$ (this is
basically the `Illusie Conjecture'---see \cite[Thm. 2.5]{J1} for a
proof).  The mapping space $\Map(U,\cI)$ may be identified with
$\Map(\Z[U],\cI)$ using adjointness, and this is just the total
complex associated to the bicomplex $(p,q)\ra \cI_q(U_p)$.  By running
the spectral sequence for the homology of this bicomplex, making use
of the fact that the $\cI_k$'s are injective sheaves ($k\geq 1$) and
$\Z[U]^\sim$ is a resolution of $\Z[X]^\sim$, one finds that the
spectral sequence collapses and the homology is just that of
$\Map(X,\cI)$.  In other words, $\Map(X,\cI) \ra \Map(U,\cI)$ is a
weak equivalence.

%%%%%%%%%%%%%%%%%%%%%%%%%%%%%%%%%%%%%%%%%%%%%%%%%%%%%%%%%%%%%%%%%%%%

\section{Other applications}
\label{se:app}

Our main application for studying hypercovers is to produce
realization functors on $\A^1$-homotopy theory \cite{DI1,Is}.  In this
section we consider a few other applications to the homotopy theory of
simplicial presheaves.

\medskip

\subsection{Change of site}\mbox{}\par
\label{subse:changesite}

Suppose that $\cC$ and $\cD$ are Grothendieck sites, and $f\colon
\cC\ra \cD$ is a functor.  The direct image functor $f_* \colon
\sPre(\cD) \ra \sPre(\cC)$ has a left adjoint $f^*$.  One is
interested in conditions on $f$ which imply that these adjoint
functors are well-behaved in relation to the homotopy theory of
simplicial presheaves.  Here is a general result which is now easy to
prove:

\begin{prop}
\label{pr:changesite}
Suppose that there is a dense set $S$ of hypercovers in $\cC$ such
that $f^*$ takes elements of $S$ to hypercovers in $\cD$.  Then the
adjoint functors $(f^*,f_*)$ give a Quillen map $U\cC_{\cL} \ra
U\cD_{\cL}$.  (Recall that a \dfn{Quillen map} is just a Quillen pair
\cite[Defn.~8.5.2]{H} regarded as a map of model categories in the
direction of the left adjoint.)
\end{prop}

In this result one {\it cannot\/} replace $U\cC_{\cL}$ by
$\sPre(\cC)_{\cL}$.  The functor $f^*$ usually does not preserve
monomorphisms, which are the cofibrations in $\sPre(\cC)_{\cL}$.

\begin{proof}
Using general facts about the universal model category $U\cC$
\cite[Prop.~2.3]{D}, the functors $(f^*,f_*)$ are a Quillen map from
$U\cC$ to $U\cD$.  If $T$ denotes the collection of hypercovers in
$\cD$, then we have assumed that $f^*$ maps $S$ into $T$.  Therefore,
by general considerations \cite[Section~5]{D} one gets a Quillen pair
between $U\cC/S$ and $U\cD/T$.  But by Theorem \ref{th:hcloc} these
localizations are just $U\cC_{\cL}$ and $U\cD_{\cL}$.
\end{proof}

Suppose that $f$ is continuous, in the sense that $\{f(U_a)\ra f(X)\}$
generates a covering sieve of $f(X)$ if $\{U_a \ra X\}$ is a covering
sieve.  It follows that $f^*$ preserves generalized covers: this is
easy for maps whose target is a representable, and the general case
can be deduced using that the target is a colimit of representables.
If one also supposes that $f^*$ preserves finite limits, then
$M_n(f^*U) \cong f^*(M_nU)$ and hence $f^*$ preserves hypercovers.
Unfortunately there are examples of interest in which $f^*$ does {\it
not\/} preserve finite limits (see \cite[Ex. 1.19, p. 103]{MV}), and
so here is a slightly different criterion which is useful:

\begin{cor}
Suppose that $\cC$ and $\cD$ are Verdier sites (see
Section~\ref{se:Verdier}).  Assume the functor $f\colon \cC\ra \cD$
preserves finite limits for diagrams of basal maps, and takes covering
families $\{U_a \ra X\}$ in $\cC$ to covering families $\{f(U_a) \ra
f(X)\}$ in $\cD$. Then $(f^*,f_*)$ give a Quillen map $U\cC_{\cL}\ra
U\cD_{\cL}$.
\end{cor}

\begin{proof}
The assumptions imply that $f^*$ preserves matching objects of basal
hypercovers (use Proposition~\ref{pr:comp1} and the material in
Section~\ref{se:Verdier}). So the condition about preserving covering
families shows that $f^*$ takes basal hypercovers in $\cC$ to basal
hypercovers in $\cD$.  Thus, Proposition~\ref{pr:changesite} applies.
\end{proof}

As an example, let $S \ra T$ be a map of schemes and consider the
base-change functor $f\colon \Sm/T \ra \Sm/S$ from the category of
smooth schemes over $T$ to the category of smooth schemes over $S$.
This functor satisfies the properties of the above proposition for any
of the standard topologies (such as Zariski, \'etale, or Nisnevich) on
$\Sm/S$ and $\Sm/T$.  So one gets a Quillen pair $U(\Sm/T)_{\cL} \ra
U(\Sm/S)_{\cL}$, by the above corollary.  Compared to the discussion
in \cite{MV}, this approach is much simpler.

\subsection{Computing homotopy classes of maps}\mbox{}\par

Given a simplicial presheaf $F$, we will use $\fibr{F}$ to denote a
fibrant replacement in $\sPre(\cC)_{\cL}$ (or in $U\cC_{\cL}$,
depending on the context).  In some sense the ultimate goal of sheaf
theory is to compute the homotopy types of the
simplicial sets $\fibr{F}(X)$.  For instance,
if $A$ is a presheaf of abelian groups and $F=K(A,n)$ is the
associated Eilenberg-Mac Lane simplicial presheaf, then
$\pi_i\fibr{F}(X)=H^{n-i}(X,\tilde{A})$ (see
Section~\ref{se:shfcohom}). If $F$ is a presheaf of chain complexes
then $\fibr{F}(X)$ computes the hypercohomology of $X$ with
coefficients in $F$, and this is where the notation $\fibr{F}$ comes
from (in the present context it goes back to \cite{Th}).

There is no known method for computing $\fibr{F}$ in general---one can
use the small object argument, but this is not very computable.  For
`nice' sites one can use the Godement resolution \cite[Prop.~3.3]{J2},
but this is also not so computable.  In this section we give analogs
of the Verdier hypercovering theorem, which show how to compute some
invariants of $\fibr{F}(X)$ using hypercovers.  It would be
interesting to construct an explicit model for the simplicial set
$\fibr{F}(X)$ using hypercovers, but we haven't been able to do this.

\medskip

We'll write $\Ho(F, G)$ for the set of weak homotopy classes of maps
from $F$ to $G$ in the homotopy category of $\sPre(\cC)_\cL$.
Likewise, $\pi(F, G)$ denotes the set $\sPre(\cC)(F, G)/\!\!\sim$,
where the equivalence relation is generated by {\em simplicial}
homotopy.

Given an object $X$ in $\cC$, let $HC_X$ denote the full subcategory
of $\sPre(\cC)$ consisting of all hypercovers of $X$.  We let $\piH_X$
denote the category with the same objects, but where 
$\piH_X(U, V)$ equals $\pi(U, V)$.

\begin{prop}
\label{pr:filt}
The category $\piH_X$ is filtered.
\end{prop}

This proposition is proved in 
\cite[Expos\'e V, 7.3.2]{SGA4} and also in
\cite[Section~8]{AM} with a slightly
different notion of hypercover (see Section \ref{se:internal}).  
We prove it again here because it is straightforward with the
techniques that we have already developed.

\begin{proof}
If $U \map X$ and $V \map X$ are both hypercovers, then so
is $U \times_X V \map X$.  Thus, we only need show that two parallel
arrows $V\dbra U$ in $\piH_X$ can be equalized.

The two maps from 
$V \map X$ to 
$U \map X$ can be assembled into the square
\[
\xymatrix{
\bd{1} \tens V \ar[r] \ar[d] & U \ar[d] \\
\del{1} \tens V \ar[r] & X              }
\]
in which the bottom map factors through $V\ra X$.  The right vertical
arrow is a local acyclic fibration by definition.  Therefore,
we apply Proposition \ref{pr:hclift} and obtain another hypercover $W
\map X$ that refines $V$, together with a diagram
\[
\xymatrix{
\bd{1} \tens W \ar[r] \ar[d] & \bd{1} \tens V \ar[r] & U \ar[d] \\
\del{1} \tens W \ar[r]\ar[urr] & \del{1} \tens V \ar[r] & X.              }
\]
The two compositions $W \map U$ are simplicially homotopic
and hence equal in $\piH_X$.
\end{proof}

The following is a generalization of the Verdier hypercovering
theorem~\cite[Expos\'e V, 7.4.1(4)]{SGA4}.  The case $K=*$ of part (b)
appeared in \cite{B}, and is cited several times in Jardine's papers
(see \cite[p. 83]{J2}, for instance).  It can be deduced from general
considerations about the category of locally fibrant simplicial
presheaves being a `category with fibrant objects'.  The
generalization to arbitrary $K$, as well as to the relative setting in
(c), doesn't seem to follow from these considerations, however.  The
case of arbitrary $K$ can be deduced from $K=*$ using
\cite[Cor. 7.5]{DI2}, but the material in Section~\ref{se:hclift}
makes it just as
easy to give a proof which handles all cases at once.

\begin{thm}
\label{th:Verdier}
Let $F$ be a locally fibrant simplicial presheaf and let $X$ belong to
$\cC$.  Let $F\ra \fibr{F}$ be a fibrant replacement for $F$ in
$\sPre(\cC)_{\cL}$.  Then
\begin{enumerate}[(a)]
\item Given a $0$-simplex $p$ of $\fibr{F}(X)$, there is a hypercover
$V \ra X$ and a map $v:V\ra F$ such that the following square commutes
up to simplicial homotopy:
\[ \xymatrix{ V \ar[r]^{v}\ar[d] & F \ar[d] \\
              X \ar[r]^{p} &\fibr{F}.}
\]
We say that `$p$ is represented by the map $v$'.  
\item Given a finite simplicial set $K$, there is an isomorphism
\[ 
\colim_{U\ra X} \pi(K\tens U,F) \map \Ho(K\tens X,F)
\]
where the colimit is taken over (the opposite category of) $\piH_X$.
\item Given $p$ and $v$ as in (a), there is an isomorphism 
\[ \pi_n(\fibr{F}(X),p) \iso 
     \colim_{U\ra V} \pi_n(\Map(U,\fibr{F}),p|_U) \iso
     \colim_{U\ra V} \pi(\del{n}/\bd{n} \tens U,F)_{v|_U}.
\]
Here $v|_U$ denotes the map $U \map V \map F$, and $p|_U$ denotes
the map $U \map V \map X \map \fibr{F}$.
The colimits are taken over the overcategory $\piH_X \downarrow V$
of hypercovers refining $V$,
and $\pi(\del{n}/\bd{n} \tens
U,F)_{v|_U}$ denotes the set of all maps $f\colon \del{n}/\bd{n}\tens U \ra
F$ such that $f\mid_{*\tens U}$ is the given map $v|_U\colon U \ra V \ra
F$, modulo simplicial homotopy relative to $*\tens U$.
\end{enumerate}
\end{thm}

\begin{proof}
Part (a) is a direct consequence of Proposition~\ref{pr:hclift2}
because $F$ and $\fibr{F}$ are both locally fibrant.

For surjectivity in (b), note that any element $\alpha$ of $\Ho(K
\tens X, F)$ is represented by an actual map $K \tens X \map
\fibr{F}$.  From Proposition \ref{pr:hclift2} again, we get a
hypercover $U \map X$ and a diagram
\[
\xymatrix{
K \tens U \ar[r]^f \ar[d] & F \ar[d] \\
K \tens X \ar[r] & \fibr{F}        }
\]
commuting up to simplicial homotopy.  
The map $f$ has image $\alpha$ in $\Ho(K \tens X, F)$.

For injectivity in (b), suppose given two maps $K \tens U \map F$ that have
the same image in $\Ho(K \tens X, F)$.  Since $K \tens U \map K \tens X$
is a local weak equivalence, this means that the two compositions
$K \tens U \map \fibr{F}$ are simplicially homotopic.  Hence we have
a diagram
\[
\xymatrix{ 
(\bd{1} \times K) \tens V \ar[r] \ar[d] & (\bd{1} \times K)
            \tens U \ar[r] & F \ar[d] \\ 
(\del{1} \times K) \tens V\ar[r]\ar@{.>}[urr] 
             & (\del{1}\times K) \tens U \ar[r] & \fibr{F} 
}
\]
for some refinement $V$ of $U$, where the lift is a
relative-homotopy-lifting.  In particular, the upper left triangle
commutes on the nose, so the two maps $K \tens U \map F$ are equal in
$\colim_{U\ra X} \pi(K\tens U,F)$.

For (c), note that the natural map $\Map(X,\fibr{F})\ra
\Map(U,\fibr{F}) $ is a weak equivalence.  So it induces an
isomorphism $\pi_n(\fibr{F}(X),p) \llra{\iso}
\pi_n(\Map(U,\fibr{F}),p|_U)$, and after taking the colimit over all
$U$ we get the first isomorphism in the theorem.

For the second isomorphism, observe that composing with $F\ra
\fibr{F}$ induces maps $\pi(\del{n}/\bd{n}\tens U,F)_{v|_U} \ra
\pi_n(\Map(U,\fibr{F}),p|_U)$.  As in the proof of part (b) above, the
fact that these maps give an isomorphism after passing to the colimit
is a direct consequence of Proposition~\ref{pr:hclift2}.
\end{proof}

Note that if $S$ is a dense set of hypercovers then the colimits in
the above results can just as well be taken over the full subcategory
of $\piH_X$ whose objects belong to $S$.

\subsection{The coconnected case}

\begin{defn}
A locally fibrant simplicial presheaf $F$ is said to be \mdfn{locally
$n$-coconnected} if it has the following property: for any $X$ in
$\cC$ and any $0$-simplex $x$ in $F(X)$, the homotopy group sheaves
$\pi_k(F,x)$ on $\cC\ovcat X$ vanish for all $k\geq n$.
\end{defn}

Using techniques from \cite{DI2},
a locally fibrant simplicial presheaf is locally $n$-coconnected
if and only if it has the local lifting property
with respect to the maps $\bd{k} \tens X \map \Delta^{k} \tens X$ for $k > n$.

Not surprisingly, for $n$-coconnected presheaves one can calculate
homotopy classes of maps by only using {\it bounded\/} hypercovers.
This is what we'll prove next.

If $n\geq 0$, let $HC_X(n)$ denote the category of bounded hypercovers
$U\ra X$ of height at most $n$ (see Definition \ref{de:bndedhc}).  Let
$\piH_X(n)$ denote the category with the same objects but with
simplicial homotopy classes of maps.  Arguments similar to
Proposition~\ref{pr:filt} show that $\piH_X(n)$ is filtered.

\begin{prop}
Suppose that $F$ is locally fibrant and locally $n$-coconnected.  Then
given a finite simplicial set $K$, there is an isomorphism
\[ \Ho(K\tens X,F) \iso \colim_{U\ra X} \pi(K\tens U,\cosk_{n} F)
\]
where the colimit is taken over the category $\piH_X(n)$.
\end{prop}

\begin{proof}
First, the map $F\ra \cosk_{n} F$ is a local weak equivalence between
locally fibrant objects.  So we can say that
\[ \Ho(K\tens X,F) \iso \Ho(K\tens X,\cosk_{n} F) \iso \colim_{U\ra X}
\pi(K\tens U,\cosk_{n} F)
\] 
where the colimit runs over the full category $\piH_X$; the second
isomorphism comes from Theorem~\ref{th:Verdier}.  We need to show that
\[ \colim_{U\in \piH_X(n)} \pi(K\tens U,\cosk_{n} F) \ra
\colim_{U\in \piH_X} \pi(K\tens U,\cosk_{n} F)
\] 
is an isomorphism.  Observe that for any simplicial set
$L$, a map $L\tens U \ra \cosk_{n} F$ factors through $\cosk_{n} (L\tens
U)$, and the map $L\tens U \ra \cosk_{n}(L\tens U)$ factors as $L\tens U
\ra L\tens \cosk_{n} U \ra \cosk_{n}(L\tens U)$.  Applying this when $L=K$
shows surjectivity because $\cosk_n U$ belongs to $HC_X(n)$, 
and from $L=K\times \del{1}$ one can deduce injectivity.
\end{proof}

\begin{prop}
\label{pr:Artin}
Let $S$ be a Noetherian scheme, and let $\cC$ be
$\Sm/S$ (or $\Sch/S$) with either the \'etale or Nisnevich
topology.  Let $X$ be an object in $\cC$ with the property that every
finite set of points is contained in an affine open.  Then every
bounded hypercover of $X$ can be refined by a \CCech complex.
\end{prop}

\begin{proof}
In the case of the \'etale topology, this is essentially the content
of \cite[Thm.~4.1]{Ar}.  Since the result is trivial for hypercovers
of height $0$, we'll suppose by induction that it works for
hypercovers of height at most $n$.  Let $U\ra X$ be a hypercover of
height $n+1$.  By Theorem~\ref{th:bhcwc}, $U$ can be refined by a
basal hypercover $U'\ra X$ (see Section~\ref{se:Verdier} below).  Let
$V=\cosk_n U'$, which is a hypercover of height at most $n$.  By
induction, there is an \'etale covering family $\{W_i \ra X\}$ such
that $\Cech{W}$ refines $V$ (where $W=\coprod W_i$).  Consider the
induced map $\Cech{W}_{n+1} \ra V_{n+1}=M_{n+1}U'$.  The map $U'_{n+1}
\ra M_{n+1}U'$ is an \'etale cover, which pulls back to an \'etale
cover $E\ra \Cech{W}_{n+1}$.  Theorem 4.1 of \cite{Ar} (applied to the
case of no geometric points) says that there is a refinement $Z$ of
$W$ such that the map $\Cech{Z}_{n+1} \ra \Cech{W}_{n+1}$ factors
through $E$.  In particular, this means that $\Cech{Z}$ refines $U'$
(and therefore $U$) up through dimension $n+1$; since the height of
$U$ is $n+1$, this means it automatically refines $U$ in all
dimensions.  This completes the proof.

For the Nisnevich topology it is essentially the same argument, only
using a revised version of \cite[Thm 4.1]{Ar}---see
\cite[Prop. 1.9, p. 99]{MV}.
\end{proof}

\begin{cor}
Let $\cC$ and $X$ be as in Proposition \ref{pr:Artin},
and suppose $F$ is a locally fibrant
simplicial presheaf which is locally $n$-coconnected.  If $K$ is a
finite simplicial set, there is an isomorphism
\[ \Ho(K\tens X,F) \iso \colim_{U\ra X} \pi(K\tens U,\cosk_{n} F)
\]
where the colimit is taken over the category $\piH_X(0)$ consisting of
the \CCech complexes.
\end{cor}

In particular, if $A$ is a presheaf of abelian groups and we take
$K=*$ and $F=K(A,n)$, then the above corollary gives the isomorphism
between \CCech cohomology and sheaf cohomology established in
\cite{Ar} and \cite[Prop. 1.9, p. 99]{MV}.

\begin{proof}
This is a direct consequence of the previous two propositions.  
The subcategory $\piH_X(0)$ is final in $\piH_X(n)$.  
\end{proof}

\begin{remark}
The above proposition and its corollary are not true for the Zariski
topology, and therefore not for the open covering topology on an
arbitrary topological space.  We repeat the example of \cite[Ex.
1.10, p. 99]{MV}: Let $X=\spec R$ be the semi-localization of $\A^2_k$ at
the points $(0,0)$ and $(0,1)$.  As a topological space $X$ has
exactly two closed points $x_1$ and $x_2$ (of codimension $2$), infinitely
many points of codimension $1$ (corresponding to the irreducible closed
curves in $\A^2$ passing through both $(0,0)$ and $(0,1)$), and a
generic point of codimension $0$.  Any open cover of $X$ can be refined by
a cover with exactly two elements: take any of the pieces containing
$x_1$ and $x_2$, respectively.

Let $U_1=X-\{x_1\}$ and $U_2=X-\{x_2\}$.  Pick two of the codimension
$1$ points $f$ and $g$ which specialize to both $x_1$ and $x_2$.  Let
$W_1=(U_1\cap U_2)-\{f\}$ and $W_2=(U_1\cap U_2)-\{g\}$.  
Let $\Omega_0=U_1\amalg U_2$ and
$\Omega_1=(U_1 \amalg U_2) \amalg W_1 \amalg W_2$ 
(the first part is degenerate).
Consider the hypercover $\cosk_1 \Omega$.
This hypercover cannot be refined by a \CCech complex.
\end{remark}

%\vfill\eject

%%%%%%%%%%%%%%%%%%%%%%%%%%%%%%%%%%%%%%%%%%%%%%%%%%%%%%%%%%%%%%%%%%%%

\section{Verdier sites}
\label{se:Verdier}

The definition of hypercover we've adopted so far in this paper has
the advantage of working for any Grothendieck site; but it is so broad
that it can sometimes be cumbersome.  One is often in the position of
having to check that something works for all hypercovers, and so it is
important to have---whenever possible---a smaller collection of
objects to deal with.  This is the subject of the present section.

To see the basic problem, look at the site of topological spaces with
the Grothendieck topology given by open covers.  Under
Definition~\ref{de:hypercover}, to give a hypercover of a space $X$
basically corresponds to giving a simplicial space $U_*$ such that
each matching map $U_n \ra M_n U$ is locally split.  This allows for
an incredible amount of freedom in what a hypercover can look like, so
much that it's very difficult to say anything concrete about it.
To make things easier, it is reasonable that one should be able to
look just at the `open hypercovers', where the maps $U_n \ra M_n U$
all have the form $\amalg_a W_a \ra M_nU$ for some open covering
$\{W_a\}$ of the target.  These are much more manageable objects.

The notion of a Verdier site---introduced in the following
definition---is just an axiomatization of the above situation.  It is
a Grothendieck site with enough extra data that one can talk about a
special kind of `basal hypercover' rather than the more general notion
we have been working with.  A Verdier site is almost just a
Grothendieck site with a basis, but we need to throw in one extra
property. 

\medskip

\begin{defn}
\label{de:Verdier}
A \dfn{Verdier site\/} is a category $\cC$ together with a given
collection of covering families $\{U_a \ra X\}$ satisfying the
properties below.  A map $U\ra X$ in $\cC$ is
\dfn{basal} if it belongs to one of these covering families.  With
this terminology, the properties can be stated as follows:
\begin{enumerate}[(i)]
\item Any single isomorphism $\{Z\ra X\}$ forms a covering family.
\item If $\{U_a \ra X\}$ is a covering family and $Y\ra X$ is a map,
then the pullbacks $Y\times_X {U_a}$ all exist, and $\{Y\times_X {U_a}
\ra Y\}$ is a covering family.
\item If $\{U_a \ra X\}$ is a covering family and one is given a
collection of covering families $\{V_{ab} \ra U_a\}$, then the
collection of compositions $\{V_{ab} \ra U_a \ra X\}$ 
is also a covering family.
\item If $U \ra X$ is a basal map then the diagonal $U\ra U\times_X U$
is also basal.
\end{enumerate}
\end{defn}

Conditions (i)--(iii) say that the collection of covering families
serves as a basis for a Grothendieck topology on $\cC$ in the usual
way.  Most of the familiar geometric Grothendieck sites satisfy the
above axioms, including topological spaces, where the covering
families are open covers, as well as the Zariski, Nisnevich, and
\'etale topologies on schemes.
The reason for not assuming that $\cC$ has all pullbacks is so that
our results apply to the Grothendieck topologies on {\it smooth\/}
schemes which are used in $\A^1$-homotopy theory \cite{MV}.

Observe that pullbacks along any basal map always exist (part (ii)),
and that any composition of basal maps is again basal (part (iii)).
It follows that if $\{V_a \ra X\}$ is a finite collection of
basal maps and $\{ U_a \map V_a \}$ is another collection of basal
maps, then the induced map $\prod_X U_a \map \prod_X V_a$ of
fibre-products is again basal.

For the following definition, note that to give a map
$f\colon\coprod_i rX_i\ra \coprod_j rY_j$ between coproducts of
representables one must choose, for every index $i$, a prescribed
value of $j$ and a map $X_i \ra Y_j$.

\begin{defn}\mbox{}\par
\label{de:basal}
\begin{enumerate}[(a)]
\item
A map $f\colon W \map Y$ in $\sPre(\cC)$ is \dfn{basal} if $W$ is a coproduct
$\coprod_i rW_i$ of representables, $Y$ is also a coproduct $\coprod_j
rY_j$ of representables, and the various maps $W_i \ra Y_j$
determining $f$ are all basal, in the sense of Definition
\ref{de:Verdier}.
\item A \dfn{basal hypercover} $U \ra X$ is a hypercover such that
the matching maps $U_n \ra M_n U$ are all basal.
\end{enumerate}
\end{defn}
The second part of this definition only makes sense if one knows that
the matching objects $M_n U$ are all coproducts of representables, but
we will see in Lemma~\ref{le:Mncoprod}
that this is the case.  First, an easy lemma:

\begin{lemma}
\label{le:basalpullback}
Let $F \ra H \la G$ be maps between coproducts of representables,
where $G\ra H$ is basal.  Then the pullback is also a coproduct of
representables, and the map from the pullback to $F$ is basal.  
\end{lemma}

\begin{proof}
Use the fact that the Yoneda embedding preserves whatever limits exist,
and that coproducts commute with fibre-products in $\sPre(\cC)$.
The necessary pullbacks in $\cC$ exist because pullbacks along basal
maps always exist.
\end{proof}

\begin{lemma}
\label{le:Mncoprod}
Let $U \map X$ be an $n$-truncated basal hypercover, and let $K$ be a
finite simplicial set of dimension at most $n$.  Then $\hom_+(K, U)$
is a coproduct of representables.  In particular, this is true for
$M_{n+1} U = \hom_+(\bd{n+1},U)$.
\end{lemma}

\begin{proof}
We proceed by induction on the dimension of $K$.  When $K$ is empty
$\hom_+(K, U)$ is just $X$, which is a representable by assumption.

Now assume that the lemma has been proven 
for simplicial sets of dimension at most $k-1$, and let 
$K$ be obtained from
a $(k-1)$-dimensional simplicial set $L$ 
by attaching finitely many $k$-simplices.  By repeating the following
argument, we may assume that only one $k$-simplex is attached.
It follows that $\hom_+ (K,U)$ is the pullback of the diagram
\[
\hom_+(\del{k},U) \map \hom_+(\bd{k},U) \leftarrow \hom_+(L, U).
\]
All three objects are coproducts of representables, the first because
$U_k$ is a coproduct of representables and the last two by the
induction hypothesis.  Since the left map above is basal (being the
matching map in a basal hypercover), Lemma \ref{le:basalpullback}
tells us that $\hom_+(K,U)$ is also a coproduct of representables.
\end{proof}

Let us return momentarily to Definition~\ref{de:basal}(b).  If $U\ra
X$ is a hypercover and $U_0 \ra X$ is basal, then the above
proposition specialized to $\bd{1}$ shows that $M_1 U$ is a coproduct
of representables.  So we may ask that $U_1 \ra M_1 U$ be basal, which
in turn forces $M_2 U$ to be a coproduct of representables.  This
shows that our definition of basal hypercover makes sense, in a
recursive sort of way.

\begin{prop}
\label{pr:basalhom}
Let $K \map L$ be any map of finite simplicial sets whose dimensions
are at most $k$, and let $U \map X$ be a $k$-truncated basal
hypercover.  Then the map $\hom_+ (L, U) \map \hom_+ (K, U)$ is basal.
\end{prop}

\begin{proof}
Consider the class $C$ of all maps of finite simplicial sets having
the property stated in the proposition.  By definition of basal
hypercovers, $C$ contains the generating cofibrations $\bd{n} \map
\del{n}$.  Cobase changes preserve $C$ by Lemmas \ref{le:hompo}(v),
\ref{le:basalpullback}, and \ref{le:Mncoprod}.  Also, finite
compositions preserve $C$ because basal maps are closed under finite
composition.  This shows that $C$ contains all inclusions of finite
simplicial sets.

In particular, $\emptyset \map \del{n}$ belongs to $C$.  This means
that $U_n \map X$ is basal for every basal hypercover $U \map X$.
Using axiom (iv) of Verdier sites one can deduce that $U_n\ra U_n
\times_X U_n$ is also basal (note that in the present context these
objects are coproducts of representables, unlike in the axiom).  In
other words, $C$ contains the codiagonal $\del{n} \amalg \del{n} \map
\del{n}$ for every $n$.

Every surjection can be built from the above codiagonals with finitely
many compositions and cobase changes.  Thus, every surjection belongs
to $C$.  But every map is a composition of a surjection with an inclusion,
so every map belongs to $C$.
\end{proof}

The proposition below is the main thing we need about basal
hypercovers.  See \cite[Lem.~8.8]{AM} for the same result without
reference to basal maps.  Unfortunately, dealing with these basal maps
definitely increases the technical complications.

\begin{thm}
\label{th:bhcwc}
In a Verdier site, any hypercover may be refined by a split, basal
hypercover.  In particular, the basal hypercovers are dense.
\end{thm}

\begin{proof}
Let $U \map X$ be any hypercover.
The fact that $U_0 \ra X$ is a generalized 
cover means there is a covering sieve
$R$ of $X$ such that every map in $R$ lifts through $U_0$.
But our Grothendieck topology was generated by a basis, so there is a
covering family $\{W_a \ra X\}$ for which every element belongs to
$R$.  Setting $V_0=\coprod_a rW_a$, we have that $V_0 \ra X$ is basal
and refines $U_0 \ra X$.

Continuing by induction, we may assume we have built a split, basal,
$n$-truncated hypercover $V$ which refines $U$ (up through dimension
$n$).  Our job is to define $V_{n+1}$.
We consider the maps
\[ \xymatrix{& U_{n+1} \ar[d] \\
          M_{n+1}V \ar[r] & M_{n+1} U,}
\]
where all the objects are coproducts of representables by 
Lemma \ref{le:Mncoprod}.
Using the same reasoning as in the first paragraph, 
there is a map $W \ra M_{n+1}V$ that is basal, that is a generalized
cover, and that fits in the upper
left corner of this diagram, i.e., it refines the pullback 
generalized cover
$U_{n+1} \times_{M_{n+1} U} M_{n+1} V \ra M_{n+1}V$.  
Set $V_{n+1}=W\amalg L_{n+1}V$.
Now $V$ is a split, $(n+1)$-truncated
hypercover; the question is whether $V_{n+1} \ra M_{n+1}V$ is basal.
Because
of the way $W$ was constructed, we need only show that
the map $L_{n+1}V \map M_{n+1} V$ is basal.

Recall from Section~\ref{se:cosk} that there is a natural map $\dgn_n V \ra
\cosk_n V$.  In dimension $n$ this is the identity map on $V_n$, and
in dimension $n+1$ it's the map $L_{n+1}V \ra M_{n+1}V$.  Picking any
degeneracy $s_i$ from level $n$ to $n+1$, we get a diagram
\[ \xymatrix{ L_{n+1} V \ar[r] & M_{n+1}V \\
              V_n \ar[u]^{s_i} \ar@{=}[r] & V_n. \ar[u]_{s_i}}
\]
Every representable summand of $L_{n+1}V$ is of the form $s_i(rU)$ for
some $i$ and some representable summand $rU$ of $V_n$, so it suffices
to show that the right-hand map $s_i: V_n \ra M_{n+1} V$ is basal.
But this degeneracy is induced by the corresponding collapse map
$\bd{n+1} \ra \del{n}$, i.e., the composition $s\colon \bd{n+1} \inc
\Delta^{n+1}\llra{s_i}\del{n}$.  In other words, $s_i$ coincides with
$\homp(\del{n},V) \ra \homp(\bd{n+1},V)$.  The fact that this is basal
follows from Proposition \ref{pr:basalhom}.
\end{proof}

\begin{remark}
\label{re:Verdier-cardinal}
Suppose there is a regular cardinal $\lambda$ with the property that
every covering family in $\cC$ has size less than $\lambda$. 
Following similar observations to
those in Proposition~\ref{pr:denseset}, the split, basal
hypercover of the above proposition
can be constructed so that in each level it has fewer than
$\lambda$ summands.  This is needed in the next section.
\end{remark}
%
%
%\vfill\eject

%%%%%%%%%%%%%%%%%%%%%%%%%%%%%%%%%%%%%%%%%%%%%%%%%%%%%%%%%%%%%%%%%%%%

\section{Internal hypercovers}
\label{se:internal}

In this final section we give a slight modification of
Theorem~\ref{th:hcloc} which is useful in applications---for instance,
it is needed in \cite{Is}.  This involves once again tweaking the
definition of hypercover in a certain way.

What sometimes happens is that the Grothendieck site $\cC$ is rich
enough that one can talk about hypercovers as elements of $s\cC$
rather than $\sPre(\cC)$, and this is usually a convenience.  For
example this is the approach taken in \cite{AM}, and it is also used
in \cite{DI1} in the context of simplicial spaces.  Handling this
involves only a slight difference from what we have done, mostly
caused by the fact that the coproduct in $\cC$ (which we will
denote by $\cup$) is not the same as the coproduct of presheaves:
i.e., $r(X\cup Y)$ is not the same as $rX \amalg rY$.

\medskip

Throughout this section we work with a Verdier site for which there
exists a regular cardinal $\lambda$ such that:
\begin{enumerate}[(1)]
\item Every covering family $\{U_i \ra X\}$ has cardinality less than
$\lambda$. 
\item Coproducts of size less than $\lambda$ exist in $\cC$.
\item If $\{X_i\}$ is a set of objects whose cardinality is less than
$\lambda$, then 
the map of presheaves $\coprod_i rX_i \ra r(\bigcup_i X_i)$ becomes an
isomorphism after sheafification.
\end{enumerate}

For example, if $\Sm/k$ denotes the category of smooth schemes of
finite type over a fixed ground field $k$, we may give it the
structure of a Verdier site by saying that the covering families are
finite collections $\{U_i \ra X\}$ such that $\coprod U_i \ra X$ is an
\'etale (or Zariski or Nisnevich) cover.  This generates the usual
Grothendieck topology, and satisfies the above properties with
$\lambda=\aleph_0$.

\begin{defn}
Given an object $X$ of $\cC$, an \dfn{internal hypercover} of $X$ is a
simplicial object $U$ in $s\cC$ which is augmented by $X$, with the
property that each matching map $U_n \ra M_n U$ is isomorphic over
$M_n U$ to a map of the form $\coprod_i V_i \ra M_n U$, for some basal
maps $\{V_i \ra M_n U\}$ which generate a covering sieve.
\end{defn}

Of course one has to worry about whether the matching object $M_n U$
exists, since the site $\cC$ need not have arbitrary limits.  But we
shall see that the condition on $U_k \ra M_k U$ for $k\leq n-1$
guarantees that $M_n U$ does in fact exist.  Even though $\cC$ is not
necessarily complete, the conclusions of Lemma~\ref{le:hompo} are
still valid when the limits $\homp(K, W)$ do exist in $\cC$.  For
example, if $\homp(L,W)$, $\homp(K,W)$, and $\homp(M,W)$ all exist,
and the pullback of 
\[ \homp(L,W) \ra \homp(K,W) \la \homp(M,W) \]
also exists, then $\homp(L\amalg_K M,W)$ exists and is isomorphic to the
above pullback.

\begin{lemma}
If $U\ra X$ is an $n$-truncated internal hypercover then the object
$\homp(K,U)$ exists whenever $K$ is a simplicial set of dimension at most $n$.
In particular, the matching object $M_{n+1} U=\homp(\bd{n+1},U)$
exists.
\end{lemma}

\begin{proof}
The proof follows the same lines as the proof of 
Lemma~\ref{le:Mncoprod}.
\end{proof}

We continue our notational convention of writing $U$ for a simplicial
object of $\cC$ and also for the simplicial presheaf that it 
represents.

\begin{thm}
The model category $\sPre(\cC)_{\cL}$ 
of simplicial presheaves may be obtained as the
localization of $\sPre(\cC)$ at the following collection of maps $\cI$:
\begin{enumerate}[(i)]
\item Maps of the form $\coprod W_i \ra (\bigcup W_i)$, for
collections $\{W_i\}$ in $\cC$ of size less than $\lambda$.
\item The maps $rU\ra rX$, for all internal hypercovers $U\ra X$.  
\end{enumerate}
\end{thm}

\begin{proof}
Let $\sPre(\cC)_{\cI}$ denote the localization we're considering.
First note that all the maps in $\cI$ are local weak equivalences.
For maps of type (ii), this is Theorem \ref{th:hcloc}.  For maps of
type (i) it follows from assumption (3) at the beginning of this
section, because every simplicial presheaf is locally weakly equivalent
to its sheafification.  So $\sPre(\cC)_{\cL}$ is a stronger
localization than $\sPre(\cC)_{\cI}$.  To see that the localizations
coincide, it will suffice to show that if $V\ra X$ is a basal
hypercover in which the number of summands in each level is smaller
than $\lambda$, then $V \ra X$ is a weak equivalence in
$\sPre(\cC)_{\cI}$.  This is by virtue of Theorem~\ref{th:hcloc},
Theorem~ \ref{th:bhcwc}, and Remark~\ref{re:Verdier-cardinal}.

Each presheaf $V_n$ may be decomposed as a coproduct of representables
in an essentially unique way: $V_n =\coprod_\alpha V_{n\alpha}$.  We
define an object $U$ of $s\cC$ by $U_n=\bigcup_\alpha V_{n\alpha}$,
and with face and degeneracy maps lifted from those in $V$.  For the
rest of the proof we will be careful to distinguish $U$ from the
simplicial presheaf $rU$.  Observe that there is a canonical map $V
\ra rU$, commuting with the augmentations down to $X$.

We claim that $U$ is an internal hypercover of $X$.  Assuming this for
the moment, relation (i) in our definition of $\sPre(\cC)_{\cI}$ shows
that $V_n \ra rU_n$ is an $\cI$-weak equivalence for each $n$.  Since
every simplicial presheaf $F$ is the homotopy colimit $\hocolim_n F_n$
(see Remark~\ref{re:hocolim}), it follows that $V \ra rU$ is also an
$\cI$-weak equivalence.  Using that $U \ra X$ is an internal
hypercover, relation (ii) gives that $rU \he X$; so one concludes
that $V \he X$ as well.

It remains only to verify that $U$ is an
internal hypercover.  First note that
$\homp(K,V)\ra \homp(K,rU)$ induces an isomorphism on sheafifications
for every finite simplicial set $K$.  When $K$ has dimension $0$ this
follows from property (3) that we assumed at the beginning of the
section.  For higher dimensional $K$ one proceeds by induction on the
number of non-degenerate simplices in $K$,
using the same pullback square from Lemma~\ref{le:Mncoprod}
and the fact that sheafification preserves finite limits.

So taking $K=\bd{n}$ we have that $M_nV \ra M_n(rU)$ induces an
isomorphism on sheafifications, and in particular is a generalized
cover.  This, together with the fact that $V_n \ra M_nV$ is a
generalized cover, shows immediately that the same is true of $rU_n \ra
M_n(rU)$.  So $rU$ is a hypercover of $X$.
 
Finally, to see that $U$ is an internal hypercover one just uses that
the Yoneda embedding preserves all limits that exist: so $M_n(rU)$ is
isomorphic to $r(M_nU)$.  
\end{proof}

%%%%%%%%%%%%%%%%%%%%%%%%%%%%%%%%%%%%%%%%%%%%%%%%%%%%%%%%%%%%%%%%%%%%

\appendix

\section{\CCech localizations}
\label{se:Cech}

This appendix is a bit of an aside from the main body of the paper.
Here we investigate how descent for \CCech complexes compares to
descent for all hypercovers.  These are not equivalent notions in
general---see Example~\ref{ex:counter}---although in some cases they
turn out to agree.  Unlike hypercover descent, \CCech descent is
often a reasonably straightforward thing to verify; so it's useful to
know how strong a notion it is.  
In this section we show that \CCech descent actually implies descent
for all {\it bounded\/} hypercovers, and we give some related results
of interest.  \CCech descent has also been explored in papers
of Simpson, for instance in \cite{HS}.  

\bigskip

If $\xi\colon F\ra G$ is a map of presheaves of sets, the \mdfn{\CCech
complex} of $\xi$ is the simplicial presheaf $\Cech{\xi}$ (often
denoted $\Cech{F}$ by abuse) given by
\[ [n]\mapsto F\times_G F\times_G \cdots
\times_G F \quad \text{($n+1$ factors)}.
\]
A simplicial presheaf $F$ is said to have \dfn{\CCech descent} if it
satisfies descent for $\Cech{U} \map X$ whenever $U\ra X$ is a generalized
cover in which $X$ is representable and $U$ is a coproduct of
representables.

Here is a short proposition we will need to use often:

\begin{prop}
\label{pr:cov-cech}
Let $\{U_a\ra X\}$ be any set of maps in $\cC$, and let $R\inc X$ be
the sieve generated by these maps.  Let $U=\coprod_a rU_a$.  
Then there is a natural map $\Cech{U} \ra R$, and this map is an
objectwise weak equivalence.
\end{prop}

\begin{proof}
If $\xi\colon K\ra L$ is any map of simplicial sets, then the \CCech
complex $\Cech \xi$ is fibrant and homotopy discrete.
This shows that the natural map $\Cech{U} \ra \pi_0 \Cech{U}$ is an
objectwise weak equivalence.  The presheaf $R$ is equal to
the presheaf $\pi_0 \Cech{U}$, i.e., $R(Y) = \pi_0
\Cech{U}(Y)$ for all $Y$ in $\cC$.
\end{proof}

Let $\ch$ denote the set of maps $\{R\inc X\}$, where $X$ runs over
all objects in $\cC$ and $R$ runs over all covering sieves (this is a
{\it set} because $\cC$ is small).  Let $\sPre(\cC)_{\ch}$ denote the
Bousfield localization of $\sPre(\cC)$ at the set $\ch$.  We'll refer to
this model category as the \dfn{\CCech localization} of $\sPre(\cC)$,
for reasons which will shortly become apparent (see
Corollary~\ref{co:gencov-cech}).

Given a covering sieve $R\inc X$, let $\Cech(R)$ denote the \CCech
complex corresponding to the cover $\coprod U_a \ra X$ where the
coproduct ranges over all maps $U_a \ra X$ in the sieve.  The above
proposition implies that $\Cech{R} \ra X$ factors through $R$, and
$\Cech{R} \ra R$ is an objectwise weak equivalence.  So localizing at
the set $\{R\inc X\}$ is equivalent to localizing at $\{\Cech{R} \ra
X\}$.  We will see in a minute that this is actually equivalent to
localizing at $\{\Cech{U}\ra X\}$, for {\it all\/} generalized covers
$U\ra X$, and so our \CCech localization is analogous to the
hypercover localization of Theorem~\ref{th:hcloc}.  The advantage of
starting with just the sieves rather than the generalized covers is
that these form a set, and so the localization automatically exists.

\begin{prop}
\label{pr:shf-cech}
Given any simplicial presheaf $F$, the map $F\ra\tilde F$ from $F$ to
its (levelwise) sheafification is a weak equivalence in
$\sPre(\cC)_{\ch}$.
\end{prop}

Unfortunately the proof of this result is somewhat involved, so we'll
postpone it until the end of the section.

\begin{cor}
\label{co:gencov-cech}
Let $F\ra G$ be any generalized cover of presheaves of sets.  Then the map
$\Cech{F} \ra G$ is a weak equivalence in $\sPre(\cC)_{\ch}$.
\end{cor}

\begin{proof}
The map $\Cech{F} \ra G$ factors as $\Cech{F} \ra \pi_0 \Cech{F} \ra
G$.  As in the proof of Proposition \ref{pr:cov-cech}, the 
first map is an objectwise weak equivalence.  The second map is a
monomorphism of presheaves, and the fact that $F\ra G$ is a
generalized cover shows that it is a local epimorphism.  Hence, the
map becomes an isomorphism upon sheafification.
Proposition~\ref{pr:shf-cech} then shows that it is a weak equivalence
in $\sPre(\cC)_{\ch}$, and so we can conclude the same for the
composite $\Cech{F}\ra G$.
\end{proof}

We now derive the connection with hypercovers.  Recall from
Definition~\ref{de:bndedhc} that a hypercover $U \map
X$ is bounded if $U = \coskp_n U$ for some $n$.

\begin{prop}
\label{pr:cech-bounded}
Given a bounded hypercover $U$ of $X$, the map
$U \ra X$ is a weak equivalence in $\sPre(\cC)_{\ch}$.
\end{prop}

The following proof was the inspiration for the proof of \cite[Lem.~4.2]{DI1}.

\begin{proof}
We proceed by induction, starting from the fact that bounded
hypercovers of height $0$ are just \CCech complexes and therefore are
handled by Corollary \ref{co:gencov-cech}.

Suppose that $U \ra X$ is a bounded hypercover of height $n+1$.
Define $V$ to be $\coskp_n U$, so $V$ is a bounded hypercover of
height at most $n$.  Therefore, we may assume by induction that $V
\map X$ is a weak equivalence in $\sPre(\cC)_{\ch}$.  The canonical
map $U \ra V$ gives a generalized cover $U_{n+1} \ra V_{n+1}$, by the
very definition of what it means for $U$ to be a hypercover.
Lemma~\ref{le:levcov} below shows that in fact $U_k \ra V_k$ is a
generalized cover for all $k$.
  
Consider the following bisimplicial object, augmented horizontally by
$V$:

\[ \xymatrix{
 V & U \ar[l] & U \times_{V} U \ar@<0.5ex>[l]\ar@<-0.5ex>[l] &
 \cdots \ar@<0.6ex>[l]\ar[l]\ar@<-0.6ex>[l] }
\]
The $k$th row is the (augmented) \CCech complex for the generalized cover
$U_k \ra V_k$.  Note that for $0\leq k \leq n$ the $k$th row
is the constant simplicial object with value $U_k$
because $U_k \ra V_k$ is the
identity.  Call this bisimplicial object (without the horizontal
augmentation) $W_{**}$.  There is an obvious map $\hocolim W_{**} \ra X$.

One may compute $\hocolim W_{**}$ by first taking the homotopy
colimit of the rows, and then taking the homotopy colimit of the
resulting simplicial object.  But in $\sPre(\cC)_{\ch}$ the homotopy
colimit of the $k$th row is just $V_k$ by Corollary~\ref{co:gencov-cech}.
Also, we have assumed by induction that
$V \he \hocolim_k V_k$ is weakly equivalent to $X$.  
So $\hocolim W_{**} \ra X$ is a weak equivalence.

Let $D$ denote the diagonal of $W_{**}$.  Standard homotopy theory
tells us that $D = \hocolim_k D_k$ is weakly equivalent to $\hocolim
W_{**}$.  We claim that $U$ is a retract over $X$ of $D$.  Note first
that one has, in complete generality, a map $U \ra D$; in dimension
$k$ it is the unique horizontal degeneracy $W_{0k} \ra W_{kk}$.

To produce a map $D \ra U$ it is enough to give $\sk_{n+1} D \ra
\sk_{n+1} U$, because $U=\coskp_{n+1} U$.  But note that $\sk_n
D = \sk_n U$, and choosing any face map $[0] \ra [n+1]$ gives a map
$D_{n+1} \ra U_{n+1}$, inducing a corresponding map $\sk_{n+1} D \ra
\sk_{n+1} U$ as desired.

It is straightforward to check that $U \ra D \ra U$ is the
identity (because $U=\coskp_{n+1} U$ one only has to check it on
$(n+1)$-skeleta), and all the maps commute with the augmentations down
to $X$.  
We have already shown that $D \map X$ is a weak equivalence in
$\sPre(\cC)_{\ch}$.
Since $U \map X$ is a retract of $D \map X$, it must also be a 
weak equivalence.
\end{proof}

\begin{lemma}
\label{le:levcov}
If $U$ is a hypercover of height $n+1$, then the map $U \ra \coskp_{n}
U$ is a generalized cover in every dimension.
\end{lemma}

\begin{proof}
First note that for $k\leq n$,
the map $U_k \ra [\coskp_n U]_k$ is the identity, so it is a 
generalized cover.
For any $k$,
the map $U_k \ra [\coskp_n U]_k$ may be rewritten as
\[\homp(\del{k},U) \ra \homp(\sk_n \del{k}, U).\]  
But $U=\coskp_{n+1} U$, so the domain may be written as
\[ \homp(\del{k},U)=\homp(\del{k},\coskp_{n+1} U)
=\homp(\sk_{n+1}\del{k},U).
\]
So we are interested in the map $\homp(\sk_{n+1}\del{k},U)\ra
\homp(\sk_n\del{k},U)$ induced by the inclusion $\sk_n\del{k}
\ra\sk_{n+1}\del{k}$.  
Recall from Lemma \ref{le:fhc} the pullback square
\[
\xymatrix{ 
\homp(\sk_{n+1} \del{k},U) \ar[d]\ar[r] 
             & \prod_X \homp(\del{n+1},U) \ar[d]\\
\homp(\sk_n \del{k},U) \ar[r] & \prod_X \homp(\bd{n+1},U). }
\]
The map $\homp(\del{n+1},U) \ra \homp(\bd{n+1},U)$ is
just the matching map $U_{n+1} \ra M_n U$, and is therefore a 
generalized cover.
So the right vertical map in the above square is a finite product of
generalized covers, which is again a generalized cover.  
Finally, we see that the left
vertical map is a pullback of a generalized cover, hence also a 
generalized cover.
\end{proof}

If $R\inc X$ is a covering sieve, let $\cI_R$ denote the full
subcategory of $\cC\ovcat X$ consisting of all maps in $R$.  Consider
the diagram $\cI_R \ra \sPre(\cC)$ sending $(U\ra X)$ to $U$.  The
colimit of this diagram is $R$, and we will write the homotopy colimit
as $\hocolim_{R} U$.  The natural map from the homotopy colimit to the
colimit gives $\hocolim_R U \ra R$, and this turns out to be an
objectwise weak equivalence by \cite[Lemma 2.7]{D}.  This fact has
nothing to do with sieves, and is true in a slightly generalized form
for arbitrary simplicial presheaves.

\begin{thm}
\label{th:allagree}
The following classes of maps give the same localization of $\sPre(\cC)$:
\begin{enumerate}[(a)]
\item The set of all covering sieves $R\inc X$.
\item The set of all maps $\hocolim_{R} U \ra X$, where $R\inc X$
is a covering sieve.
\item The class of all hypercovers of height $0$, i.e., the 
\CCech complexes $\Cech{U} \ra X$.
\item The class of all bounded hypercovers $U \ra X$.
\item The class of maps $F \ra \tilde F$ from simplicial presheaves to 
their sheafifications.
\end{enumerate}

\noindent
If the topology on $\cC$ is given by a basis of covering families,
then one can also add
\begin{enumerate}[(a')]
\item The set of all covering sieves $R_U\inc X$ where $R_U$ is the
sieve generated by the covering family $\{U_a \ra X\}$. 
\end{enumerate}
\end{thm}

It may seem surprising that the localization in (e) does not give the
usual notion of local weak equivalence, but only the weaker \CCech
version.  Example~\ref{ex:counter} shows that it really
is weaker.  Also, note that the above theorem could just as well have
been stated for $U\cC$ rather than $\sPre(\cC)$---the proofs are
essentially the same.

\begin{proof}
The fact that $\hocolim_R U \ra R$ is an objectwise weak equivalence
immediately shows that the localizations in (a) and (b) are the same.
And we have seen in Proposition~\ref{pr:cov-cech} that the
localizations in (a) and (c) are the same.

The localization in (d) is {\it a priori\/} stronger than that in (a);
Proposition~\ref{pr:cech-bounded} shows that the two localizations 
actually agree.  Likewise, the localization in (e) is stronger than
the one in (a), because $R\inc X$ becomes an isomorphism upon
sheafification.  Proposition~\ref{pr:shf-cech} shows that they agree.

Finally, if our topology is given by a basis of covering families then
the proof that the localizations in (a') and (e) coincide follows the
proof of Proposition~\ref{pr:shf-cech} more or less verbatim.
\end{proof}

It would be interesting to know more about $\sPre(\cC)_{\ch}$, for
instance to have an explicit characterization of the weak equivalences.  
Perhaps this wouldn't be so useful, since 
the chief interest in
$\sPre(\cC)_{\ch}$ is that it is sometimes a more convenient version
of $\sPre(\cC)_{\cL}$ (see Example \ref{ex:ch=L}).

\begin{cor}
\label{co:bhcdesc}
Let $F$ be a simplicial presheaf.  Then $F$ satisfies descent for all
\CCech complexes if and only if it satisfies descent for all bounded
hypercovers.
\end{cor}

\begin{proof}
Let $F'$ be a fibrant replacement for $F$ in $\sPre(\cC)$.  The
statement of the corollary for $F'$ requires that $F'$ be local with
respect to the \CCech complexes $\Cech{V} \ra X$ if and only if it is
local with respect to the bounded hypercovers $U \ra X$.  This is true
by the above proposition (parts (c) and (d)).  But of course $F$ has
descent for a certain class of objects precisely when $F'$ has descent
for that same class, because $F\ra F'$ is an objectwise weak
equivalence.
\end{proof}

\begin{cor}
\label{co:cechbasis}
Suppose the Grothendieck topology on $\cC$ is given by a basis of
covering families.  Then a simplicial presheaf $F$ satisfies \CCech
descent if and only if it satisfies descent for all the \CCech
complexes $\Cech{U}\ra X$ in which $X$ is a representable and
$U=\coprod_a U_a$ for some covering family $\{U_a \ra X\}$ in the
basis.
\end{cor}

\begin{proof}
This is similar to the proof of Corollary~\ref{co:bhcdesc}, using
Theorem~\ref{th:allagree} (parts (a') and (c)) and
Proposition~\ref{pr:cov-cech}.
\end{proof}

The following result can be useful for verifying hypercover descent.
We deduce it from our general analysis of hypercovers, but the
statement also follows from results of \cite{HS} (see especially
Proposition 6.1) and should be credited to that paper.

\begin{cor}
Let $F$ be an objectwise-fibrant simplicial presheaf with the property
that $F(X)$ has no homotopy in dimension $n$ or higher, for every $X$ in
$\cC$.  Then $F$ satisfies descent for all hypercovers if and only if
it satisfies descent for all \CCech complexes.
\end{cor}

\begin{proof}
First we need to consider the localization $U\cC/S$, where $S$ is the
set of maps $\{\bd{n+1}\tens X \ra \del{n+1}\tens X| X\in \cC\}$.  It
is easy to check that the fibrant objects in $U\cC/S$ are the
simplicial presheaves $G$ such that each $G(X)$ is fibrant and has no
homotopy above dimension $n-1$.  Given an objectwise-fibrant
simplicial presheaf $G$, one can construct the fibrant replacement $L_S G$
via the small object argument applied to the maps in $S$.  By thinking
about this, one sees that the maps of simplicial sets $G(X) \ra L_S
G(X)$ are isomorphisms up through simplicial dimension $n$.  So $L_S
G(X)$ has the same homotopy groups as $G(X)$ up through dimension
$n-1$, but no homotopy groups in higher dimensions.  Even if $G$ is
not objectwise-fibrant, $L_SG$ is objectwise weakly equivalent to
$L_S(\Ex^\infty G)$, and so it is still true that $G(X)$ and $L_SG(X)$
have the same $(n-1)$-type for all $X\in \cC$.

General localization theory says that a map $G\ra H$ is an
$S$-equivalence if and only if $L_S G \ra L_S H$ is an objectwise
equivalence, and so this is the same as saying that $G(X)\ra H(X)$
induces isomorphisms on all homotopy groups up through dimension
$n-1$, for every $X$.  In particular, the map $G\ra \cosk_n G$ is
an $S$-equivalence.

Now consider the localization $U\cC/T$, where $T$ is the union
of $S$ and the set of all covering sieves $R\inc X$.  A simplicial
presheaf $F$ is fibrant in $U\cC/T$ precisely when it is 
objectwise-fibrant, has \CCech descent, and each $F(X)$ has no homotopy in
dimension $n$ or higher.

Suppose that $F$ is as in the statement of the corollary, and that $F$
satisfies descent for all \CCech complexes.  Then we know that $F$ is
fibrant in $U\cC/T$.  But note that $U\ra \cosk_n U$ is necessarily a
$T$-equivalence (because it is an $S$-equivalence).  Yet $\cosk_n U$
is a bounded hypercover of $X$, and hence $\cosk_n U \ra X$ is a
$T$-equivalence as well (using the $U\cC$ version of
Proposition~\ref{pr:cech-bounded}).  Hence $U\ra X$ is also a
$T$-equivalence.  Since $F$ is fibrant in $U\cC/T$ and $X$ is
cofibrant in $U\cC/T$, the morphism $\Map(X,F)\ra \Map(U',F)$ is a
weak equivalence, where $U'$ is a cofibrant replacement for $U$ in
$U\cC$.  Thus $F$ satisfies descent for $U\ra X$ by Lemma
\ref{le:descent}(ii).
\end{proof}

Here is an example showing that the \CCech localization can be
strictly weaker than the localization at all hypercovers.  In other
words, we exhibit a simplicial presheaf which has descent for all
\CCech complexes but does not have descent for all hypercovers.  The
example is a slight modification of one suggested to us by Carlos
Simpson.

\begin{example}
\label{ex:counter}
Let $X=X_0$ be the open interval $(0,1)$.  
Now let $U_0=(0,\frac{2}{3})$, $V_0=(\frac{1}{3},1)$,
and $X_1=U_0\cap V_0$.  Note that $X_1 \iso X$, and let
$U_1=(\frac{1}{3},\frac{5}{9})$ , $V_1=(\frac{4}{9},\frac{2}{3})$, and
$X_2=U_1\cap V_1$.  Again one has $X_2 \iso X$, and we define $U_2$,
$V_2$, and $X_3$ in the expected way.  Continue.  Our site $\cC$
consists of the spaces $\{X_i,U_i,V_i \mid i \geq 0\}$ with the
inclusions between them, and equipped with the usual notion of open
cover.  The category $\cC$ is depicted as
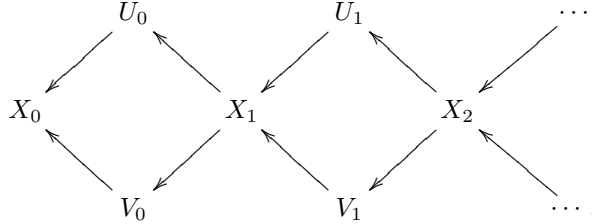
\[
\xymatrix{
     & U_0 \ar[dl] & &  U_1 \ar[dl] & &  \cdots \ar[dl] \\
X_0 & & X_1 \ar[ul] \ar[dl] & & X_2 \ar[ul] \ar[dl] \\
     & V_0 \ar[ul] & &  V_1 \ar[ul] & &  \cdots. \ar[ul] \\
}
\]
Define a presheaf of topological spaces on our site in the following way:
\[ F(X_0)=\emptyset,
\quad F(U_n)=D^n_{+},
\quad F(V_n)=D^n_{-},
\quad \text{and} \quad F(X_{n+1})=S^n \ \ (n\geq 0).
\]
Here $D^n_{+}$ and $D^n_{-}$ denote the upper and lower hemispheres of
$S^n$.  The restriction maps $F(U_n) \ra F(X_{n+1})$ and $F(V_n) \ra
F(X_{n+1})$ are the inclusions of the hemispheres in $S^n$, while the
maps $F(X_n) \ra F(U_n)$ and $F(X_n) \ra F(V_n)$ are the inclusions of
the boundaries of the hemispheres.  Define the simplicial presheaf $G$
by $G(W)=\Z \Sing F(W)$---that is, $G(W)$ is the result of applying
the free abelian group functor to the singular complex of $F(W)$.
Using the Dold-Kan correspondence, one can regard $G$ as a presheaf of
chain complexes; then $G(W)$ is the usual complex for computing the
singular homology of $F(W)$.

Now $G$ has \CCech descent: this can be checked using
Corollary~\ref{co:cechbasis}, and so the main point is that for every
$n$ the square
\[\xymatrix{ G(U_n\cup V_n) \ar[r] \ar[d] & G(U_n) \ar[d] \\
              G(V_n) \ar[r] & G(U_n \cap V_n)} 
\]
is a homotopy pullback.  
% Here we really need the free group functor.  The homotopy pullback
% of $F(U_n) \ra F(X_{n+1}) \la F(V_n)$ is $\Omega S^n$, which is not
% $S^{n-1}$.
On the other hand, we will construct a hypercover for which $G$ does
not have descent.  The combinatorics of this construction are slightly
complicated, but the idea is this: Start with $\Omega_0 = U_0 \amalg
V_0$, then consider $\cosk_0 \Omega$ except replace each
non-degenerate occurrence of $X_1$ with $U_1 \amalg V_1$.  Next take
$\cosk_1 \Omega$, replace each non-degenerate occurrence of $X_2$ with
$U_2 \amalg V_2$, and continue.  This gives the hypercover $\Omega\ra
X$.

Now we will be more precise.  Let $P_n$ be the category of all
nonempty subsets of $\{0,1,\ldots,n\}$, with inclusions.  Note that
the objects of $P_n$ can be identified with the sub-simplices of
$\del{n}$, and that $[n]\mapsto P_n$ forms a simplicial category.  Let
$S_n$ denote the set of all functors $J\colon P_n^{op} \ra \cC$ with
the following properties:

\begin{enumerate}[(1)]
\item All the values of $J$ belong to $\{U_i,V_i\mid i\geq 0\}$.
\item Given a subset $\sigma=\{i_0,\ldots,i_k\}$ in $P_n$, if
$\bigcap_j J(\{i_0,\ldots,\hat{i_j},\ldots,i_k\}) = U_m$ (resp. $V_m$)
then $J(\sigma)=U_m$ (resp. $V_m$).
% In other words, if $m$ is the largest index that occurs on the faces and
% if $U_m$ occurs on a face but $V_m$ does not, then
% $J(\sigma)$ is also $U_m$.
\item If $\bigcap_j J(\{i_0,\ldots,\hat{i_j},\ldots,i_k\})= X_m$ then
$J(\sigma)$ is either $U_m$ or $V_m$ (and this includes the case $k=0$).
% In other words, if $m$ is the largest index that occurs on the faces and
% if both $U_m$ and $V_m$ occur on faces, then
% $J(\sigma)$ is either $U_{m+1}$ or $V_{m+1}$.
\end{enumerate}

Let $\Omega$ denote the simplicial presheaf defined by
\[ 
   \Omega_n := \coprod_{J\in S_n} J(\{0,1,\ldots,n\}),
\]
with simplicial structure induced by that of $P$.  Intuitively, each
summand of $\Omega_n$ corresponds to an $n$-simplex together with a
certain labelling of its simplices given by $J$: the labelling is such
that smaller simplices are labelled by larger opens, and such that the
above properties are satisfied.  The reader is encouraged to work out
what these properties say for small values of $n$, and to verify that
$\Omega$ is a hypercover of $X_0$ (use Proposition~\ref{pr:comp1}).
Check that $\Omega_0=U_0\amalg V_0$ and $\Omega_1=U_0 \amalg (U_1
\amalg V_1) \amalg (U_1\amalg V_1) \amalg V_0$.  

We claim that $\holim_n G(\Omega_n)$ is not connected, whereas
$G(X)=0$.  To calculate $\pi_0(\holim_n G(\Omega_n))$ we can just work
in the category of chain complexes.  The cosimplicial object $[n]
\mapsto G(\Omega_n)$ corresponds to a double complex, and we are
trying to compute the $0$th homology of the total complex (the one
called $\Tot^\Pi$ in \cite{W}, rather than $\Tot^\oplus$).  But
observe that each $G(\Omega_n)$ has homology only in dimension $0$
because each $G(U_i)$ and $G(V_i)$ is contractible.  Therefore, the
$E_1$-term of the spectral sequence for the homology of the bicomplex
is concentrated in a line.  Thus, the bicomplex's $0$th homology is
the kernel of $d_0 - d_1\colon H_0 G(\Omega_0) \map H_0 G(\Omega_1)$,
which is $\Z$.  This completes the verification that $G$ does not
satisfy descent for the hypercover $\Omega$.
\end{example}

\begin{example}
\label{ex:ch=L}
Sometimes the localizations $U\cC_{\ch}$ and $U\cC_{\cL}$ {\it do\/}
coincide.  Let $S$ be a Noetherian scheme of finite dimension and let
$\cC$ be the site $\Sm/S$ with either the Zariski or Nisnevich
topology (one can also take $\Sch/S$ here).  For these sites the
localizations $U\cC_{\ch}$ and $U\cC_{\cL}$ agree.  For the
Zariski topology this is a direct consequence of the `Brown-Gersten
Theorem', which identifies the fibrant objects in $U\cC_{\cL}$ with
the objectwise-fibrant simplicial presheaves satisfying \CCech descent for
all two-fold Zariski covers $\{U,V\ra X\}$.  This is essentially
proven in \cite{BG}, although one has to translate their proof into
our more modern setting.  See also \cite[Lem. 4.1, 4.3]{Bl}.

For the Nisnevich topology we have to explain a little more.  Given an
elementary distinguished square $\{U\inc X,p\colon V\ra X\}$
\cite[Def. 1.3, p. 96]{MV}, let $P(U,V)$ denote the simplicial presheaf
which has $U\amalg V$ in dimension $0$, $U\amalg p^{-1}(U) \amalg V$
in dimension $1$, and is degenerate in higher dimensions.  The
Brown-Gersten Theorem in this context is \cite[Lem. 1.6, p. 98]{MV};
together with
\cite[Lem. 4.3]{Bl}, it implies that $U\cC_{\cL}$ is the
localization of $U\cC$ at the maps $P(U,V) \ra X$ for all elementary
distinguished squares.  
We already know that $U\cC_{\cL}$ is a stronger localization than
$U\cC_{\ch}$, so we just need to show that the maps
$P(U,V) \map X$ are weak equivalences in $U\cC_{\ch}$.

To see that $P(U,V)\ra X$ is a weak equivalence in $U\cC_{\ch}$, first
note that for any $Z$ the simplicial set $P(U,V)(Z)$ has
non-degenerate simplices only in dimensions $0$ and $1$.  Each
component is a star, centered at a $0$-simplex corresponding to a map
$Z\ra U$ (because every map $Z\ra V$ can be an endpoint of at most one
$1$-simplex).  Therefore each component is contractible, so $P(U,V)
\ra \pi_0P(U,V)$ is an objectwise weak equivalence.  The codomain is
just the presheaf $U\amalg_{p^{-1}U} V$, so we are reduced to showing
that the map $U\amalg_{p^{-1}U} V\ra X$ is a weak equivalence in
$U\cC_{\ch}$.  By the $U\cC$ version of Theorem~\ref{th:allagree}(e)
it suffices to show that this map induces an isomorphism on
sheafifications, and this is routine (use Nisnevich stalks, or look at
\cite[Lem. 1.6, p. 98]{MV}).
\end{example}

\subsection{A leftover proof}
\label{se:leftover}
The final goal of this section is to give the proof of
Proposition~\ref{pr:shf-cech}: if $F$ is a simplicial presheaf we need
to show that $F\ra \tilde F$ is a weak equivalence in
$\sPre(\cC)_{\ch}$.  In fact it will suffice to do this when $F$ is a
discrete simplicial presheaf, since a simplicial presheaf $F$ can be
recovered as a homotopy colimit of the discrete presheaves $F_n$
(Remark \ref{re:hocolim}).  Unfortunately, even to prove our claim for
discrete simplicial presheaves seems to require a wrestling match with
the small object argument.

So from now on $F$ is just a presheaf of sets.  We introduce two
constructions: First, $\cA F$ is the presheaf defined by $\cA F(X) =
F(X)/\sim$, where two sections $s$ and $t$ are equivalent if there is
a covering sieve $R\inc X$ such that $s|_U=t|_U$ for every $U\ra X$ in
$R$.  
Secondly, $\cB F$ is defined to be the pushout
\[ \xymatrix{ \coprod R \ar[r] \ar[d] & F \ar[d] \\
              \coprod X \ar[r] & \cB F
}
\]
where the coproduct is indexed over all objects $X$ in $\cC$, all
covering sieves $R\inc rX$, and all maps $R \ra F$.  One may check
that $\cA\cB F$ is what is usually denoted $F^+$, and so $\cA\cB\cA\cB
F$ is the sheafification $\tilde{F}$.

We will show that the maps $F \ra \cB F$ and $F\ra \cA F$ are \CCech
weak equivalences, for any presheaf $F$.  The first claim is very
easy: Since $\coprod R \ra \coprod X$ is an acyclic cofibration in
$\sPre(\cC)_{\ch}$, its cobase change $F\ra \cB F$ is also an acyclic 
cofibration.  Unfortunately the
second claim is much more difficult.  The idea is to build up a
\CCech weak equivalence $F \ra L_\infty F$ by brute force, in such a
way that there is an objectwise weak equivalence $L_\infty F \ra \cA
F$.

Given a covering sieve $R\inc X$, let $J^n R$ be the pushout
\[ \xymatrix{ \bd{n} \tens R \ar[d] \ar[r] & \bd{n}\tens X \ar[d] \\
              \del{n} \tens R \ar[r] & J^n R.}
\]
The natural map $J^n R \ra \del{n} \tens X$ is a \CCech acyclic
cofibration
because $\sPre(\cC)_{\ch}$ is a simplicial model category.
Note that a map $J^{n} R \ra G$ is the
same as a map $\bd{n} \ra G(X)$ together with a compatible
family of extensions $\del{n} \ra G(U)$ for all maps $U\ra X$ in $R$.

Let $L_0 F=F$, and let $L_{n+1}F$ be obtained from $L_n F$ as the
pushout
\begin{myequation}
\label{eq:Ln-po}
 \xymatrix{ \coprod J^{n+1} R \ar[r] \ar@{ >->}[d]_\sim 
                              & L_n F \ar@{ >->}[d]^\sim \\
               \coprod \del{n+1} \tens X \ar[r] & L_{n+1} F
}
\end{myequation}
where the coproducts run over all $X$ in $\cC$, all covering sieves
$R\inc X$, and all maps $J^{n+1} R \ra L_n F$.  Let $L_\infty F$ be
the colimit of the chain $L_0 F \ra L_1 F \ra L_2 F \ra \cdots$.
Since each $L_n F \ra L_{n+1} F$ is a \CCech acyclic cofibration, the
composite $F \ra L_\infty F$ is also a \CCech acyclic cofibration.

To get a feel for what's happening here, let's look just at $L_1F$.  A
map $J^1 R \ra F$ corresponds to giving a map $\bd{1}\ra F(X)$
together with a compatible family of extensions $\del{1}\ra F(U)$ for
all $U\ra X$ in the sieve.  Since $F$ is discrete, this means that we
are giving two elements of $F(X)$ which {\it agree\/} when restricted
to pieces of the sieve.  When we form the pushout $\del{1}\tens X \la
J^1 R \ra F$ we are adding a $1$-simplex into $F(X)$ which
identifies these elements in $\pi_0$.  So it follows that $\pi_0
L_1F(Y)=\cA F(Y)$, for all $Y$.

When we pass from $L_nF$ to $L_{n+1}F$ something similar is
happening---we will see it boils down to killing off all the higher
homotopy, in the end because the objects $F(X)$ were all discrete and
so had no higher homotopy to begin with.  So the goal is to show that
each $L_\infty F(X)$ is fibrant and homotopy discrete, and that $\pi_0
L_\infty F(X)=\cA F(X)$.  This will imply that the natural map
$L_\infty F \ra \pi_0 L_\infty F$ is an objectwise weak equivalence,
and the target is identified with $\cA F$.  We will then have $F\we
L_\infty F \we \cA F$ in $\sPre(\cC)_{\ch}$.

The argument will proceed by establishing the following properties:
Given any $Y$ in $\cC$,
\begin{enumerate}[(i)]
\item The map of simplicial sets $L_n F(Y) \ra L_{n+1} F(Y)$ is an
isomorphism on $n$-skeleta.
\item The simplicial set $L_n F(Y)$ has dimension at most $n$, i.e.,
it is degenerate in degrees greater than $n$.
\item Given any $n$-simplex $\sigma$ in $L_nF(Y)$, there is a covering sieve
$R\inc Y$ such that $\sigma|_{U}$ is in the image of $L_{n-1}F(U) \ra
L_n F(U)$ for every $U\ra Y$ in $R$.  In particular, $\sigma|_U$ is
a degenerate $n$-simplex.
\item  Given any $n$-simplex $\sigma$ in $L_nF(Y)$, there is a
covering sieve $R\inc Y$ 
such that $\sigma|_{U}$ is in the image of $F(U) \ra
L_n F(U)$ for every $U\ra Y$ in $R$.
\item For $n\geq 2$, any map $\bd{n} \ra L_{n-1} F(Y)$ extends to a
map $\del{n} \ra L_n F(Y)$.
\item Any map $\Lambda^{2,k} \ra L_1 F(Y)$ extends to $\bd{2} \ra L_2
F(Y)$.
\end{enumerate}

Granting these for the moment, let us show they imply the desired
result.  To show that $L_\infty F(X)$ is fibrant and homotopy
discrete, it is enough to verify that it has the extension property
with respect to the maps $\bd{n} \ra \del{n}$ ($n\geq 2$) and the maps
$\Lambda^{2,k} \ra \del{2}$.  These are easy consequences of parts
(i), (v), and (vi).  Also, part (i) tells us that $\pi_0 L_1 F \ra
\pi_0 L_\infty F$ is an isomorphism, and we have already remarked that
$\pi_0 L_1 F \cong \cA F$.  This finishes the proof, granting the
statements outlined above.

Claim (i) follows from the fact that $J^{n+1} R(Y) \ra (\del{n+1}\tens
X)(Y)$ is an isomorphism on $n$-skeleta.  Part (ii) follows from an
induction, using that $(\del{n}\tens X)(Y)$ has dimension $n$ and that
$F(Y)$ has dimension $0$ (since we assumed that $F$ is a presheaf of
sets).  Part (iii) is a straightforward analysis of diagram
(\ref{eq:Ln-po}), and (iv) follows from (iii) by induction.  We will
show that (v) is a consequence of (iv), and a similar argument proves
(vi).

Suppose we have a map $\sigma\colon\bd{n} \ra L_{n-1}F(Y)$.  By (iv),
for each face $d_i\sigma$ there is a covering sieve $R_i \inc Y$
such that $d_i\sigma|_{U}$ is in the image of $F(U)\ra L_{n-1}F(U)$,
for every $U\ra Y$ in $R_i$.  There is of course a covering sieve $R$
which refines all the $R_i$.  So for each $U\ra Y$ in $R$ and each
$i$, there is $(n-1)$-simplex $\alpha_{U,i}$ in $F(U)$ which maps to
$(d_i\sigma)|_{U}$.

Now, it is not clear that as $i$ varies the $(n-1)$-simplices
$\alpha_{U,i}$ fit together to give a map $\alpha_U\colon \bd{n} \ra
F(U)$.  However, we know they fit together in $L_{n-1}F(U)$, and the
map $F(U) \ra L_{n-1}F(U)$ is a cofibration of simplicial sets, hence
a monomorphism.  So the $\alpha_{U,i}$ must fit together in $F(U)$ as
well.

Secondly, it is not immediately clear that the $\alpha_U$ patch
together over the covering sieve $R$: that is, we must check that
given maps $U\ra V \ra X$ where $V\ra X$ is in $R$, then $\alpha_U$
coincides with the restriction of $\alpha_V$ to $U$.  Again, this
follows from the fact that everything patches together in $L_{n-1}F$
and the fact that $F \ra L_{n-1}F$ is an objectwise cofibration.

So we have constructed a map $\alpha\colon \bd{n}\tens R \ra F$ such
that the composite map $\bd{n}\tens R \ra F \ra L_{n-1}F$ coincides with
$\sigma|_R$.  Now we use the fact that $F$ is a {\it discrete} simplicial
presheaf, from which it follows that $\alpha$ can be extended to a map
$\bar{\alpha}\colon \del{n} \tens R \ra F$.  Composing this with $F
\ra L_{n-1}F$ and patching with $\sigma$ gives a map $J^n R \ra
L_{n-1}F$, and this extends to $\del{n}\tens Y$ once we pass to
$L_n F$.  The upshot is that we've shown $\sigma$ extends to $\del{n}$
under the map $L_{n-1}F(Y)\ra L_{n} F(Y)$.

%\vfill\eject

%%%%%%%%%%%%%%%%%%%%%%%%%%%%%%%%%%%%%%%%%%%%%%%%%%%%%%%%%%%%%%%%%%%%

\bibliographystyle{amsalpha}

\end{document}